\documentclass[11pt]{article}
\usepackage[english]{babel}
\usepackage{amsmath,amssymb,bbm,pslatex,xcolor,graphics}

\usepackage{fancyhdr}

\catcode`\<=\active \def<{
\fontencoding{T1}\selectfont\symbol{60}\fontencoding{\encodingdefault}}
\catcode`\>=\active \def>{
\fontencoding{T1}\selectfont\symbol{62}\fontencoding{\encodingdefault}}
\catcode`\|=\active \def|{
\fontencoding{T1}\selectfont\symbol{124}\fontencoding{\encodingdefault}}
\newcommand{\tmrsub}[1]{\ensuremath{_{\textrm{#1}}}}

\newcommand{\gott}{\ensuremath{\mathfrak{T}}}

\newcommand{\tptwo}{\ensuremath{T^{\ast} ( \varnothing )_{p_{2}}}}
\newcommand{\pe}{\ensuremath{P^{\ensuremath{\operatorname{eJ}}}_{e}}}
\newcommand{\tmt}[1]{\ensuremath{\mathbbm{T}^{#1}}}
\newcommand{\p}[1]{\ensuremath{P^{\ensuremath{\operatorname{eJ}}}_{#1}}}
\newcommand{\zx}[1]{\ensuremath{^{#1} \zeta}}
\newcommand{\sx}[1]{\ensuremath{^{#1} \Sigma}}

\newcommand{\go}{G{\"o}del}

\newcommand{\phiabp}[2]{\ensuremath{\Phi_{\beta +1}^{(  \alpha )}}}

\newcommand{\con}{{\smallfrown}}

\newcommand{\ptwo}{{{\em II\/}}}
\newcommand{\pone}{{{\em I\/}}}

\newcommand{\nod}{{\noindent}}

\newcommand{\bu}{\ensuremath{\bullet}}
\newcommand{\nat}{\ensuremath{\mathbbm{N}}}

\newcommand{\bai}{\textsuperscript{\ensuremath{\omega}}\ensuremath{\omega}}

\newcommand{\cant}{2\textsuperscript{\ensuremath{\mathbbm{N}}}}
\newcommand{\can}{\textsuperscript{\ensuremath{\omega}}2}
\newcommand{\rest}{{\upharpoonright}}
\newcommand{\emp}{{\varnothing}}
\newcommand{\pa}[1]{\ensuremath{\langle #1 \rangle}}

\newcommand{\game}{\ensuremath{\Game}}

\newcommand{\da}{{\downarrow}}
\newcommand{\ua}{{\uparrow}}
\newcommand{\la}{{\langle}}
\newcommand{\ra}{{\rangle}}

\newcommand{\back}{{\backslash}}
\newcommand{\power}{P}

\newcommand{\ie}{{\itshape{i.e.}}{\hspace{0.25em}}}
\newcommand{\eg}{{\itshape{e.g.}}{\hspace{0.25em}}}
\newcommand{\etc}{{\itshape{etc.}}{\hspace{0.25em}}}

\newcommand{\all}{{\forall}}

\newcommand{\equi}{{\longleftrightarrow}}
\newcommand{\ex}{{\exists}}

\newcommand{\sset}{{\subseteq}}
\newcommand{\lr}{{\leftrightarrow}}
\newcommand{\imp}{{\longrightarrow}}

\newcommand{\fin}[1]{\ensuremath{[  ]^{< \omega}}}

\newcommand{\qed}{Q.E.D.}

\newcommand{\pf}{{\noindent}{\textbf{Proof: }}}
\newcommand{\dfs}{\ensuremath{=_{\ensuremath{\operatorname{df}}}}}

\newcommand{\oddpagetext}[1]{\newcommand{\pageoddheader}{{\small }}}
\newcommand{\evenpagetext}[1]{\newcommand{\pageevenheader}{{\small }}}

\catcode`\<=\active \def<{
\fontencoding{T1}\selectfont\symbol{60}\fontencoding{\encodingdefault}}
\catcode`\>=\active \def>{
\fontencoding{T1}\selectfont\symbol{62}\fontencoding{\encodingdefault}}
\catcode`\|=\active \def|{
\fontencoding{T1}\selectfont\symbol{124}\fontencoding{\encodingdefault}}
\newcommand{\asterisk}{*}
\newcommand{\nobracket}{}
\newcommand{\nocomma}{}

\newcommand{\tmem}[1]{{\em #1\/}}
\newcommand{\tmop}[1]{\ensuremath{\operatorname{#1}}}
\newtheorem{theorem}{Theorem}[section]

\newtheorem{corollary}[theorem]{Corollary}
\newtheorem{definition}[theorem]{Definition}

\newtheorem{lemma}[theorem]{Lemma}
\newtheorem{remark}[theorem]{Remark}
\newtheorem{example}[theorem]{Example}
\newtheorem{note}[theorem]{Note}

\def\pone{$I$}
\begin{document}
\title{$G_{\delta \sigma}$-games and generalized computation }
\author{ P.D. Welch\\
\,\,\, School of Mathematics, University of Bristol, \\Bristol, BS8 1TW, England}
\maketitle
\begin{abstract}
  We show the equivalence between the existence of winning strategies for
  $G_{\delta \sigma}$ (also called $\Sigma^{0}_{3}$) games in Cantor or Baire
  space, and the existence of functions generalized-recursive in a higher
  type-2 functional. (Such recursions are associated with certain transfinite
  computational models.)
  
  We show, {{\em inter alia,\/}} that the set of indices of convergent
  recursions in this sense is a complete $\Game \Sigma_{3}^{0}$ set: as
  paraphrase, the listing of those games at this level that are won by player
  \pone, essentially has the same information as the `halting problem' for
  this notion of recursion.
  
  Moreover the strategies for the first player in such games are recursive in
  this sense. We thereby establish the ordinal length of monotone $\Game
  \Sigma^{0}_{3}$-inductive operators, and characterise the first ordinal
  where such strategies are to be found in the constructible hierarchy. In
  summary:
  
  {\textbf{Theorem}} (a) \ {{\em The following sets are recursively
  isomorphic.\/}}
  
  (i)  {{\em The complete ittm-semi-recursive-in-$\tmop{eJ}$ set,\/}}
  $H^{\ensuremath{\operatorname{eJ}}}$;
  
  (ii) {{\em the $\Sigma_{1}$-theory of\/}} $( L_{\eta_{0}} , \in ) ,$ {{\em
  where $\eta_{0}$ is the closure ordinal of $\Game \Sigma_{3}^{0}$-monotone
  inductions\/}};
  
  (iii) {{\em the complete $\Game \Sigma_{3}^{0}$ set of integers.\/}}
  
  \ \ \ \ (b) {{\em The ittm-recursive-in-$\tmop{eJ}$ sets of integers are
  precisely those of\/}} $L_{\eta_{0}}$.
\end{abstract}
\section{Introduction}

The attempt to prove the determinacy of two person perfect information games
(and the consequences of the existence of such winning strategies) has a long
and fruitful history, starting with work of Banach and Mazur and continuing to
the present. The work in the paper {\cite{We11}} was initially motivated by
trying to see how the $\Pi^{1}_{3}$-theory of {\itshape{arithmetical
quasi-inductive definitions}} fits in with other subsystems of second order
number theory, in particular with the determinacy of $\Sigma^{0}_{3}$-sets.
There it was shown, {{\em inter alia\/}}, that AQI's - which were known to be
formally equivalent with the most basic form of generalized computation to be
introduced below - are not strong enough to compute strategies for
$\Sigma^{0}_{3}$-games. What had been left open was a more precise discussion
of the location of those strategies. We continue that discussion here. To give
this research a context we shall also mention the results previously known in
this area.

The argument in {\cite{We12}} explicitly extracts what was undeclared in the
proof, a criterion for where exactly the strategies appear in the {\go}
constructible $L_{\alpha}$ hierarchy. Whilst we have had this result for some
while, the characterisation is somewhat unusual in that it is expressed in
terms of the potential for such $L_{\alpha}$ to have certain kinds of
ill-founded elementary end extensions, and is not so perspicuous. We had
conjectured that certain kinds of illfounded-computation trees (defined by
Lubarsky) should also characterize this ordinal. This we have verified, but
now see that there is a bigger picture that connects the generalized recursion
theory of the late 50's and early 60's of Kleene (v.{\cite{Kl59}}) of higher
types with the determinacy of games at this level. To be clearer the
connection is between the existence of winning strategies and the {{\em
generalization\/}} of Kleene which is associated with a transfinite
computational model of the so-called Infinite Time Turing machines of Hamkins
and Kidder {\cite{HL}}. Kleene in {\cite{Kl59}} developed an equational
calculus, itself evolving out of his analysis of the {\go}-Herbrand General
Recursive Functions (on integers) fom the 1930's, but now enlarged for dealing
with recursion in objects of finite type. (The set of natural numbers we
denote by $\omega$ and they are of type 0; $f:a \rightarrow \omega$ is of type
$k+1$ if $a$ is of type $k$.) A particular type-2 functional was that of the
{{\em ordinary jump\/}} $J$, where
$$
J ( e, \vec{m} , \vec{x} )  = \left\{ \begin{array}{ll} 1 & \mbox { if } \{ e \} ( \vec{m} , \vec{x} ) \da \mbox{ (meaning {{\em
converges\/}}, or {{\em is defined\/}}) } \\
 0 & \mbox{ otherwise.} 
 \end{array}\right.
$$

Here $\vec{m}$ is a string of integers, and $\vec{x}$ a vector of functions
$f: \omega \rightarrow \omega$ (thus a vector of objects of type 1) and $\{ e
\}$ a usual index of a recursive function. The function under discussion is
$\{ e \}$ which is given by a natural number index coding its formation. In
this formalism the index set

$$H^{J} ( e ) \lr \{ e \}^{J} ( e ) \da$$

{\nod}is a complete semi-recursive (in $J$) set of integers, and Kleene showed
that this is in turn {\nod}a complete $\Pi^{1}_{1}$ set of integers. \ Further
he showed that the {{\em $J$-recursive sets of integers\/}}, {\ie} those sets
$R $ for which
$$R ( n ) \lr \{ e
\}^{J} ( n ) \da  1 \, \wedge\,\neg R ( n ) \lr \{ e \}^{J} ( n ) \da  0$$

{\nod}for some index $e$, are precisely the hyperarithmetic ones.

Recall that a set $X \sset \omega   ( \bai )$ is said to be in $\game \Gamma$
for some (adequate) pointclass $\Gamma$ on the integers (Baire space), if
there is a set $Y \sset \omega \times \bai   \, ( \bai \times \bai )$ so that
$X=\{x \mid  ${{\em Player I has a winning strategy in $G (Y_{x} ,  ^{<
\omega} \omega ) \}$ \/}}where $Y_{x} = \{y  \mid \pa{x,y} \in Y\}$. Roughly
speaking, if one has a recursive listing of the $\Gamma$ sets of reals, (say
from some universal $\Gamma$ set): $A_{0} ,A_{1}  , \ldots  ,A_{n} , \ldots$ ,
then a {{\em complete\/}} $\game \Gamma$ set of integers, gives those $n$ for
which {\pone} has a winning strategy in $G ( A_{n} ; ^{< \omega} \omega )$.

We have the following theorem connecting this with determinacy of open games:

\begin{theorem}
  (Moschovakis {{\em {\cite{Mosch71}}\/}}, Svenonius {{\em {\cite{Sve65}}\/}})
  The complete $\game \Sigma^{0}_{1}$ set of integers is a complete
  $\Pi^{1}_{1}$ set of integers.
\end{theorem}

Hence by Kleene's results just alluded to:

\begin{corollary}
  The complete $\game \Sigma^{0}_{1}$ set of integers is recursively
  isomorphic to $H^{J}$, a complete $J$-semi-decidable set of integers.
\end{corollary}

Moreover:

\begin{theorem}
  (Blass {{\em {\cite{Bla72}}\/}}) Any $\Sigma^{0}_{1}$-game for which the
  open player, that is {\pone}, has a winning strategy, has a hyperarithmetic
  winning strategy.
\end{theorem}

\begin{corollary}
  Any $\Sigma^{0}_{1}$-game for which player {\pone} has a winning strategy,
  has a $J$-recursive strategy.
\end{corollary}

We seek to raise these ideas to the level of $\Sigma^{0}_{3}$. Kleene also
gave an equivalent account of recursion in objects of finite type using as an
alternative the Turing model enhanced with oracle calls to a higher type
functional, see {\cite{Kl62b}},{\cite{Kl62a}}; the account here is motivated
in spirit by that approach. Instead of using an equational calculus we shall
couch this in terms of {{\em infinite time Turing machines\/}} -(ittm's)
computations recursive in a certain operator $\ensuremath{\operatorname{eJ}}$
in place of $J$. Indeed there is already a version of this kind of computation
in the literature. In {\cite{Lu10}} Lubarsky defines the notion of a {{\em
`feedback'-\/}}{{\em ittm\/}}{{\em  machine\/}}, where a Hamkins-Kidder ittm
may call upon a sub-computation handled by another such machine, and pass an
index and an element of Cantor space to it as a parameter. The information
passed back is as to whether the computation with the given index acting on
the given parameter {{\em halts\/}} or not (which it may do after a
transfinite number of steps, in contradistinction to the standard Turing
machine). \ This is thus in the spirit of the jump $J$ defined above. A
convergent {{\em feedback-ittm computation\/}} can then be conceived as a
wellfounded tree of halting sub-computations. A {{\em divergent
computation\/}} (``freezing'' in Lubarsky's terminology) is one which descends
down an ill-founded path.

Rather than define recursions involving what would be the generalization of
$J$ above to {{\em halting\/}} ittm-computations, we use an{{\em  eventual
jump\/}} operator $\ensuremath{\operatorname{eJ}}$. The ittm's have an
arguably more fundamental behaviour than `halting' or `non-halting': they may
eventually have some settled output on their output tape without formally
entering a halting state (the Read/Write head may be meandering up and down
the tape, perhaps fiddling with the Scratch or Input tape, but leaving the
output alone, in some fixed loop without formally halting). This `eventual' or
`settled' behaviour fits in with the $\Sigma_{2}$ definable liminf rules of
its operation. We thus define:
$$\ensuremath{\operatorname{eJ}} ( e, \vec{m} , \vec{x} )  = 
\left\{ \begin{array}{ll}
1 & \mbox{ if }\{ e \} (
\vec{m} , \vec{x} ) \mid \mbox{ (denoting {{\em converges to a settled ouput\/}})} \\
 0 & \mbox{  otherwise. }
\end{array} \right.
$$

Here $\{ e \}$ is now an index of a standard ittm-computable function, say
given by some usual finite programme $P_{e} ( \vec{m} , \vec{x} )$. We then
consider ittm-computations recursive in $\ensuremath{\operatorname{eJ}}$, for
which we would now use the notation $\{ e \}^{\ensuremath{\operatorname{eJ}}}$
to denote the $e$'th such function recursive in
$\ensuremath{\operatorname{eJ}}$. Here a {{\em query\/}} instruction or state
is included as part of the machine's language. For this notion we find a level
of the $L$ hierarchy $L_{\alpha_{0}}$ to provide an analogy with the above.

\begin{theorem}
  The complete $\game \Sigma^{0}_{3}$ set of integers is recursively
  isomorphic to $H^{\ensuremath{\operatorname{eJ}}}$, the complete
  $\ensuremath{\operatorname{eJ}}$-semi-decidable set of integers.
\end{theorem}

Thus to paraphrase, the listing of those games that are won by {\pone},
essentially has the same information as the `halting problem' for this notion
of recursion. We feel this is interesting as it demonstrates that two, {{\em
prima facie\/}} very different, notions are in fact intimately connected.
Define $\tau_{0}$ as the supremum of the convergence times of
$\ensuremath{\operatorname{eJ}}$-recursive computations.

Corresponding to the result on $\Pi^{1}_{1}$ we have:

\begin{theorem}
  The complete $\game \Sigma^{0}_{3}$ set of integers is a complete
  $\Sigma^{L_{\tau_{0}}}_{1}$ truth set.
\end{theorem}

(Recall that the complete $\Pi^{1}_{1}$ set is also the
$\Sigma^{L_{\omega_{1}^{\ensuremath{\operatorname{ck}}}}}_{1}$ truth set.)
$_{}$ Moreover

\begin{theorem}
  Any $\Sigma^{0}_{3}$-game for which the player {\pone} has a winning
  strategy, has an $\ensuremath{\operatorname{eJ}}$-recursive winning
  strategy.
\end{theorem}

Corresponding to the result on hyperarithmetic strategies we have:

\begin{corollary}
  Any $\Sigma^{0}_{3}$-game for which player {\pone} has a winning strategy,
  has a winning strategy in $L_{\tau_{0}}$.
\end{corollary}

We assume the reader has familiarity both with the constructible hierarchy of
G{\"o}del - for which see Devlin {\cite{De}}. For the basic notions of
descriptive set theory including the elementary theory of Gale-Stewart games,
see Moschovakis {\cite{Mosch3}}. Our notation is standard. Some of the results
here relate to sub-systems of second order number, or analysis, and the basic
theory of this is exposited in Simpson's monograph {\cite{Si99}}. For models
of admissible set theory, also called ``Kripke-Platek set theory'' or ``KP''
see Barwise {\cite{Bar}}. By ``KPI'' we mean the theory KP augmented by the
axiom that every set is an element of some admissible set.

In the language of generalized recursion theory, the pointclass $\Game
\Sigma^{0}_{3}$ of sets of integers cannot be the 1-envelope of a normal
type-2 function, by results of Harrington, Kechris, and Simpson (see
{\cite{HaKe75}}). (A ``1-envelope'' is the set of relations on $\omega$
recursive in the type-2 functional.) What we are showing here is that the
complete set of integers in $\Game \Sigma^{0}_{3}$ is however (recursively
isomorphic to) the complete set which is ittm-semi-recursive in
$\ensuremath{\operatorname{eJ}}$ - the eventual jump type-2 functional. It is
the ``ittm-1-envelope'' of $\ensuremath{\operatorname{eJ}}$. Section 3
contains some facts related to ittm-computations, and an exposition, and sets
the scene with some basic results of our
ittm-recursions-in-$\ensuremath{\operatorname{eJ}}$.

We answer a further question of Lubarsky concerning Freezing-ITTM's at
Corollary \ref{Bob}.\\

{\nod}Acknowledgements: We should like to warmly thank Bob Lubarsky for
illuminating explanations of his paper {\cite{Lu10}}, discussions on the
conjecture mentioned in the second paragraph, and helpful comments on an
earlier draft of this paper.

\section{}

We first repeat the extraction from our earlier paper \cite{We12} of a
criterion for the constructible rank of $\Sigma^{0}_{3}$ games' strategies.
(Note that we take our games as defined in $L$ and using constructible, indeed
an initial recursive, game trees; the existence of a winning strategy for a
particular $\Sigma^{0}_{3}$ (indeed arithmetic or Borel) game is a
$\Sigma^{1}_{2}$ assertion about the countable tree $T$ and the payoff set. As
$T \in L$ the truth of such an assertion has the same truth value in the
universe of sets or in $L$. We thus expect to find such strategies in $L$
(since Davis in {\cite{Da}} proved such strategies exist in the universe $V$
of sets). \ But where are they?

\begin{definition}
  A pair of ordinals $( \mu , \nu $) is a $\Sigma_{2}${{\em -extendible
  pair\/}}, if $L_{\mu} \prec_{\Sigma_{2}} L_{\nu}$ and moreover $\nu$ is the
  least such with this property. We say $\mu$ is $\Sigma_{2}${{\em
  -extendible\/}} if there exists $\nu$ with $( \mu , \nu $) a
  $\Sigma_{2}${{\em -extendible pair\/}}. By relativisation, a pair of
  ordinals $( \mu , \nu $) is an $x$-$\Sigma_{2}$-{{\em extendible pair\/}},
  and $\mu$ is $x$-$\Sigma_{2}$-{{\em extendible\/}}, if $L_{\mu} [ x ]
  \prec_{\Sigma_{2}} L_{\nu} [ x ]$.
\end{definition}

Indeed all the above ideas relativi{\textbf{}}se normally to real parameters
$x \in \cant$, and we thus have $\lambda ( x ) , \zeta ( x ) , \Sigma ( x )$
{\etc}, with the latter two forming the least $x$-$\Sigma_{2}$-extendible
pair.

\begin{definition}
  Let an $m$-{{\em depth $\Sigma_{2}$-nesting\/}} of an ordinal $\alpha$ be a
  sequence $( \zeta_{n} , \sigma_{n} )_{n<m}$ with (i) For $0 \leq n<m$:
  $\zeta_{n-1} \leq \zeta_{n} < \alpha < \sigma_{n} < \sigma_{n-1}$ ; \ \ (ii)
  $L_{\zeta_{n}} \prec_{\Sigma_{2}} L_{\sigma_{n}}$. \ We write $d ( \alpha )
  \geq m$. If $\alpha$ is not nested we set $d ( \alpha ) =0$.
\end{definition}

We shall want to consider non-standard admissible models $(M,E)$ of
$\ensuremath{\operatorname{KP}}$ together with some other properties. We let
$\ensuremath{\operatorname{WFP}} (M)$ be the wellfounded part of the model. By
the so-called `Truncation Lemma' it is well known ({{\em v.\/}} {\cite{Bar}})
that this well founded part must also be an admissible set. Usually for us the
model will also be a countable one of ``$V=L$''. \ Let $M$ be such and let
$\alpha =\ensuremath{\operatorname{On}} \cap \ensuremath{\operatorname{WFP}}
(M)$. By the above $\alpha$ is thus an `admissible ordinal', {\ie}
$L_{\alpha}$ will also be a $\ensuremath{\operatorname{KP}}$ model. An
`$\omega$-depth' nesting cannot exist by the wellfoundedness of the ordinals.
However an ill founded model $M$ when viewed from the outside may have
infinite descending chains of `$M$-ordinals' in its ill founded part. These
considerations motivate the following definition.

\begin{definition}
  An {{\em infinite depth $\Sigma_{2}$-nesting\/}} {{\em of \/}}$\alpha$ {{\em
  based on $M$\/}} is a sequence $( \zeta_{n} ,s_{n} )_{n< \omega}$ with :
  
  (i) $\zeta_{n-1} \leq \zeta_{n} < \alpha \subset s_{n} \subset s_{n-1}$ ;
  (ii) $s_{n} \in \ensuremath{\operatorname{On}}^{M} ; $ \ (iii)
  $(L_{\zeta_{n}} \prec_{\Sigma_{2}} L_{s_{n}} )^{M}$.
\end{definition}
Thus the $s_{n}$ form an infinite descending $E$-chain through the illfounded
part of the model $M$. In {\cite{We11}} we devised a game whereby one player
produced an $\omega$-model of a theory and the other player tried to find such
infinite descending chains through $M$'s ordinals. In this paper we shall
switch the roles of the players, and have Player {{\em II\/}} produce the
model and Player $I$ attempt to find the chain. (This is just to orientate the
game as then $\Sigma^{0}_{3}$.)

In order for there to exist a non-standard model with an infinite depth
nesting (of the ordinal of its wellfounded part) then the wellfounded part
will already be a relatively long countable initial segment of $L$ (it is easy
to see that if $\zeta = \sup_{n} \zeta_{n}$ then already $L_{\zeta} \models
\Sigma_{1}$-Separation).

\begin{example}
  (i) Let $\delta$ be least so that $L_{\delta} \models
  \Sigma_{2}$-Separation, and let $(M,E)$ be an admissible non-wellfounded end
  extension of $L_{\delta}$ with $L_{\delta}$ as its wellfounded part. Then
  there is an infinite depth nesting of $\delta$ based on $M$.
  
  (ii) By refining considerations of the last example, let $\gamma_{0}$ be
  least such that there is $\gamma_{1} > \gamma_{0}$ with $L_{\gamma_{0}}
  \prec_{\Sigma_{2}} L_{\gamma_{1}} \models \ensuremath{\operatorname{KP}}$.
  Then again there is an infinite depth nesting of $\gamma_{1}$ based on some
  illfounded end extension $M$ of $L_{\gamma_{1}}$.
\end{example}

Both of the above can be established by standard Barwise Compactness
arguments. However both these $\delta$ and $\gamma_{0}$ we shall see are
greater than the ordinal $\beta_{0}$ defined from this notion of nesting as
follows.

\begin{definition}
  Let $\beta_{0}$ be the least ordinal $\beta$ so that $L_{\beta}$ has an
  admissible end-extension $(M,E)$ based on which there exists an infinite
  depth $\Sigma_{2}$-nesting of $\beta$. 
\end{definition}

\begin{definition}
  Let $\gamma_{0}$ be the least ordinal so that for any game $G (A, T)$ with
  $A \in \Sigma^{0}_{3}$, $T \in L_{\gamma_{0}}$ a game tree, then there is a
  winning strategy for a player definable over $L_{\gamma_{0}}$.
\end{definition}

The following then pins down the location of winning strategies for games at
this level played in, {\eg} recursive trees.

\begin{theorem}
  \label{gamma=beta}$\gamma_{0} = \beta_{0}$. Moreover, any
  $\Sigma^{0}_{3}$-game for a tree $T$, with a strategy for Player I, has such
  a strategy an element of $L_{\beta_{0}}$. Any $\Pi^{0}_{3}$-game for such a
  tree has a strategy which is definable over $L_{\beta_{0}}$.
\end{theorem}

\begin{definition}
  Let $\eta_{0}$ be the closure ordinal of monotone $\game
  \Sigma_{3}^{0}$-operators.
\end{definition}

This ordinal will be less than $\beta_{0}$.

\begin{theorem}
  \label{2.8}(a) \ The following sets are recursively isomorphic.
  
  (i)  The complete ittm-semi-recursive-in-$\ensuremath{\operatorname{eJ}}$
  set, $H^{\ensuremath{\operatorname{eJ}}}$;
  
  (ii) the $\Sigma_{1}$-theory of$( L_{\eta_{0}} , \in )$;
  
  (iii) the complete $\Game \Sigma_{3}^{0}$ set of integers.
  
  {\nod}(b) The ittm-recursive-in-$\ensuremath{\operatorname{eJ}}$ sets of
  integers are precisely those of $L_{\eta_{0}}$. 
\end{theorem}

\begin{definition}
  Let $\tau_{0}$ be the supremum of convergence ordinals of well-founded
  computations, arising from infinite time Turing machine computations on
  integers which are ittm-recursive (in a generalized sense of Kleene {{\em et
  al.\/}}) in the Type-2 eventual jump functional
  $\ensuremath{\operatorname{eJ}}$.
\end{definition}

\begin{theorem}
  $\label{sigma=tau} \eta_{0} = \tau_{0}$.
\end{theorem}

Remark: (i) The proof reveals more about the $L$-least strategies for
$\Sigma^{0}_{3}$-games: those for player {\pone}, in fact can be found within
a strictly bounded initial segment of $\beta_{0}$: they will occur in
$L_{\eta_{0}}$.

(ii) The existence of all the above ordinals, and $\beta$-models of the above
theories can be proven in the subsystem of analysis
$\Pi^{1}_{3}$-$\ensuremath{\operatorname{CA}}_{0}$, but not in
$\Delta^{1}_{3}$-$\ensuremath{\operatorname{CA}}_{0}$ (or even some
strengthenings of the latter). See {\cite{We11}}.

\subsection{The location of strategies for $\Sigma^{0}_{3}$-games }

{\nod}{\pf} {\textbf{of Theorem \ref{gamma=beta}}} \ We look at the
construction of the proof of Theorem 5 \ of {\cite{We11}} in particular that
of Lemma 3. There we used an assumption that there is a triple of ordinals
$\gamma_{0} < \gamma_{1} < \gamma_{2}$ with (a) $L_{\gamma_{0}}
\prec_{\Sigma_{2}} L_{\gamma_{1}}$ and (b) $L_{\gamma_{0}} \prec_{\Sigma_{1}}
L_{\gamma_{2}}$ and (c) $\gamma_{2}$ was the second admissible ordinal beyond
$\gamma_{1}$. \ One assumed that $I$ did not have a winning strategy in $G
(A;T)$. The Lemma 3 there ran as follows:

\begin{lemma}
  \label{2.12} Let $B \subseteq A \subseteq \lceil T \rceil$ with $B \in
  \Pi^{0}_{2}$. If $(G (A;T)$ is not a win for $I)_{L_{\gamma_{0}}}$, then
  there is a quasi-strategy $T^{\ast} \in L_{\gamma_{0}}$ for {\ptwo} with the
  following properties:
  
  (i) $\lceil T^{\ast} \rceil \cap B =  \varnothing \text{ }$
  
  (ii) $ ( G (A;T^{\ast} )$ is not a win for $I )_{L_{\gamma_{0}}}$.
\end{lemma}

The format of the lemma's proof involved showing that the
$\Sigma^{L_{\gamma_{0}}}_{2}$ notion of `goodness' embodied in (i) and (ii)
held for $\emp$. \ To do this involved defining goodness in general. We first
define $T'$ as {\ptwo}'s {{\em non-losing quasi-strategy\/}} for $G (A;T)$
(the set of positions $p \in T$ so that {\pone} does not have a winning
strategy in $G ( A;T_{p} )$); this is $\Pi_{1}$ definable over
$L_{\gamma_{0}}$ as the latter is a model KPI; in particular if we use the
notation

\begin{definition}
  $S^{1}_{\gamma} \dfs \{\delta < \gamma \mid  L_{\delta} \prec_{\Sigma_{1}}
  L_{\gamma} \}$.
\end{definition}

\begin{definition}
  For $n \leq \omega$, let $T^{n}_{\delta}$ denote the $\Sigma_{n}$-theory of
  $L_{\delta}$.
\end{definition}

{\noindent}then `` $p \in T'$ '' is $\Pi^{L_{\zeta_{0}}}_{1}   $, where
$\zeta_{0} \dfs \min  S^{1}_{\gamma_{0}} \back \rho_{L} (T)$. More generally
we define:

A$ p  \in  T'  $is {{\em good\/}} if there is a quasi-strategy $T^{\ast}  $
for {{\em II\/}} in $T'_{p}  $ so that the following hold:

(i) $\lceil T^{\ast} \rceil   \cap  B =  \varnothing$; \

(ii) $G (A;T^{\ast} )$ is not a win for $I$.

Here $T'_{p}$ is the subtree of $T'$ below the node $p$. \ The point of
requiring that the pair $( \gamma_{0} , \gamma_{1} )$ have the
$\Sigma_{2}$-reflecting property of (a) above, is that the class $H$ of good
$p$'s of $L_{\gamma_{1}}$ is the same as that of $L_{\gamma_{0}}$ and so is a
set in $L_{\gamma_{1}}$ as it is thus definable over $L_{\gamma_{0}}$ by a
$\Sigma_{2} ( \{ T' \} )$ definition. \ The overall argument is a proof by
contradiction, where we assume that $\emp$ is in fact not good, and proceeds
to construct a strategy $\sigma$ for Player $I$ in the game $G (A;T' )$, which
is definable over $L_{\gamma_{1}}$, and is apparently winning in
$L_{\gamma_{2}}$. (The requirement (c) that $\gamma_{2}$ be a couple of
admissibles beyond $\gamma_{1}$ was only to allow for the strategy $\sigma$ to
be seen to be truly winning by going to the next admissible set, and verifying
that there are no winning runs of play for {\ptwo}.) The contradiction arises
since $T'$ - which was defined as the subtree of $ T$ of {\ptwo}'s non-losing
positions - is concluded still to be the same subtree of non-losing positions
in $L_{\gamma_{2}}$. Being a non-losing position, $p$ say, for {\ptwo} is a
$\Pi_{1}$ property of $p$. This carries up from $L_{\gamma_{0}}$ to
$L_{\gamma_{2}}$ as $L_{\gamma_{0}} \prec_{\Sigma_{1}} L_{\gamma_{2}}$, and
this is the reason for the requirement (b): we want $T'$ to survive beyond
$L_{\gamma_{1}}$ for our argument to work. (This idea is important for the
arguments in Section 4, so let us refer to it as `{{\em the survival
argument\/}}'.) \ There is then no winning strategy for $I$ in $G (A;T' )$
definable over $L_{\gamma_{1}}$, contradicting the reasoning that $\sigma$ is
such.

This proves the Lemma: $L_{\gamma_{1}}$ sees there is $T^{\ast}$ a subtree of
$T'$ witnessing that $\emp$ is good. The existence of such a subtree is a
$\Sigma_{2} ( \{ T' \} )$-sentence, and then again this reflects down to
$L_{\gamma_{0}}$. We thus have such a $T^{\ast}$ in $L_{\gamma_{0}}$.

The Theorem is proven by repeated applications of the Lemma, by using the
argument for each $\Pi^{0}_{2}$ set $B_{n}$ in turn where $A= \bigcup_{n}
B_{n}$ and refining the trees using this procession from a tree to a subtree
$T^{\ast}$. We thus repeat the argument with $T^{\ast}$ replacing$ T$. \
Because $T^{\ast} \in L_{\gamma_{0}}$ we have the same constellation of this
triple of ordinals $\gamma_{i}$ above the constructible rank of $T^{\ast}$,
and can do this.

However we can get away with less. The definition of the subtree of non-losing
positions of {\ptwo} now this time in the new $T^{\ast}$ can be considered as
taking place $\Pi_{1}$ over $L_{\delta_{0}}$ where $\delta_{0}$ is the least
element of $S^{1}_{\gamma_{0}}$ with $T^{\ast} \in L_{\delta_{0}}$. \ To get
our contradiction we actually use that $L_{\delta_{0}} \prec_{\Sigma_{1}}
L_{\gamma_{2}}$ ; we do not need that $L_{\gamma_{0}} \prec_{\Sigma_{1}}
L_{\gamma_{2}}$. \ Notice that our argument that $T^{\ast}  $exists is
non-constructive: we simply say that the $\Sigma_{2}$-sentence of its
existence reflects to $L_{\gamma_{0}}$: we do not have any control over its
constructible rank below $\gamma_{0}$. \ Moreover any sufficiently large
$\gamma'$ greater than $\gamma_{1}$ would do for the upper ordinal, as long as
it is a couple of admissibles larger than $\gamma_{1}$. Thus we could apply
the Lemma repeatedly for different $B_{n}$ if we have a guarantee that
whenever a $T_{n}^{\ast}$-like subtree is defined there exists a $\zeta_{n}
\in S^{1}_{\gamma_{0}}$ and a suitable upper ordinal $\gamma_{n} > \gamma_{1}$
with $T_{n}^{\ast} \in L_{\zeta_{n}} \prec_{\Sigma_{1}} L_{\gamma_{n}}$ . Of
course if there are arbitrarily large $\zeta_{n}$ below $\gamma_{0}$ with this
extendability property, then this is tantamount to $L_{\gamma_{0}}
\prec_{\Sigma_{1}} L_{\gamma'}$ for some suitable $\gamma'$, and this shows
why our original constellation of $\gamma_{i}$ provides a sufficient
condition.

Actually as the final paragraph of the Theorem 5 there shows, we are doing
slightly more than this: we are, each time, applying the Lemma infinitely
often to each possible subtree of $T^{\ast}$ below some node $p_{2}$ of it
which is of length $2$, to define our strategy $\tau$ applied to moves of
length $3$. \ We then move on to the next $\Pi^{0}_{2}$ set. Although we are
applying the Lemma infinitely many times for each such $p_{2}$, and thus
infinitely many new $\Sigma_{2}$-sentences, or trees, have to be instantiated,
we had that $L_{\gamma_{0}}$ is a $\Sigma_{2}$-admissible set, and as the
class of such $p_{2}$ is just a set of $L_{\gamma_{0}}$, \
$\Sigma_{2}$-admissibility works for us to find a bound for the ranks of the
witnessing trees, as some $\delta < \gamma_{0}$. \ We thus can claim that our
final $\tau$ is an element of $L_{\gamma_{0}}$ even after $\omega$-many
iterations of this process.

$( \beta_{0} \geq \gamma_{0} )$ We argue for this. Let $(M,E)$ be a
non-standard model of $\ensuremath{\operatorname{KP}}$ with an infinite
nesting $( \zeta_{n} ,s_{n} )$ about $\beta_{0}$ as described. Note that
$S^{1}_{\beta_{0}}$ must be unbounded in $\beta_{0}$ (so that $L_{\beta_{0}}
\models \Sigma_{1}$-Separation), and each $\zeta_{n}$ is a limit point of
$S^{1}_{\beta_{0}}$. We do not assume that $\beta_{0}$ is
$\Sigma_{2}$-admissible (which in fact it is not as the proof shows). Let $T
\in L_{\beta_{0}}$ be a game tree. By omitting finitely much of the outer
nesting we assume $T \in L_{\zeta_{0}}$. \ We assume that Player $I$ has no
winning strategy for $G (A;T)$ in $L_{\beta_{0}}$ (for otherwise we are done).
Note that in $M$ we have that $L_{s_{0}}$ also has no winning strategy for
this game (otherwise the existence of such would reflect into$ L_{\beta_{0}}$.
We show that $I I$ has a winning strategy definable over $L_{\beta_{0}}$. Let
$A= \bigcup B_{n}$ with each $B_{n} \in \Pi^{0}_{2}$. For $n=0$ we apply the
argument of the Lemma using the pair $( \zeta_{1} ,s_{1} )$ in the role of $(
\gamma_{0} , \gamma_{1} )$ from before, with $( \zeta_{0} ,s_{0} )$ in the
role of $( \delta_{0} , \gamma_{2} )$ described above, {\ie} we use only that
$T_{} \in L_{\zeta_{0}}$ and that $L_{\zeta_{0}} \prec_{\Sigma_{1}}
L_{s_{0}}$.

The Lemma then asserts the existence of a quasi-strategy for {\itshape{II}}
definable using the pair $( \zeta_{1} ,s_{1} )$: $T^{*} ( \varnothing )$. By
$\Sigma_{2}$-reflection the $L$-least such lies in $L_{\zeta_{1}}$, and we
shall assume that $T^{*} ( \varnothing )$ refers to it.\\

{{\em Claim: For any pair $( \zeta_{n} ,s_{n} ) \text{ for } n \geq 1$ the
same tree $T^{\asterisk} ( \varnothing )$ would have resulted using this
pair.\/}}

Proof: \ Note that we can define such a tree like $T^{*} ( \varnothing )$
using such pairs, since for all of them we have that $( \zeta_{0} ,s_{0} )
\supset ( \zeta_{1} ,s_{1} ) \supset ( \zeta_{m} ,s_{m} )$ for $m>1$. As
$T^{*} ( \varnothing ) \in L_{\zeta_{1}}$ and satisfies a $\Sigma_{2}$
defining condition there, and since we also have $\zeta_{1} \in
S^{1}_{\zeta_{m}}$, it thus satisfies the same $\Sigma_{2}$ condition in
$L_{\zeta_{m}}$. \ {\hspace*{\fill}}{\qed} {{\em Claim\/}}\\

For any position $p_{1} \in T$ with $\ensuremath{\operatorname{lh}} (p_{1} )
=1$, let $\tau (p_{1} )  $be some arbitrary but fixed move in $T' (
\varnothing )$, this now {\itshape{II}}'s non-losing quasi-strategy for the
game $G (A,T^{*} ( \varnothing ))$ as defined in $L_{\zeta_{2}}$. \ The
relation ``$p \in T' ( \varnothing )$'' is $\Pi_{1}^{L_{\zeta_{2}}} (\{T^{*} (
\varnothing )\})$ or equivalently $\Pi_{1}^{L_{\zeta_{1}}} (\{T^{*} (
\varnothing )\})$, or indeed $\Pi_{1}^{L_{\delta}} (\{T^{*} ( \varnothing
)\})$ where $\delta$ is least in $S^{1}_{\zeta_{1}}$ above $\rho_{L} (T^{*} (
\varnothing ))$. \ Hence ``$y=T' ( \varnothing )$'' {{\em $\in
\Delta^{L_{\delta}}_{2} (\{$\/}}$T^{*} ( \varnothing ) \})$ and thus $T' (
\varnothing )$ also lies in $L_{\zeta_{1}}$. For definiteness we let $\tau
(p_{1} )$ be the numerically least move.

For any play, $p_{2}$ say, of length 2 consistent with the above definition
of $\tau$ so far, we apply the lemma again with $B=A_{1}$ replacing $B=A_{0}$
and with $(T^{*} ( \varnothing ))_{p_{2}}$ replacing $T$. We use the nested
pair $( \zeta_{2} ,s_{2} )$ to define quasi-strategies for {\itshape{II,}}
call them $T^{*} (p_{2} )$, one for each of the countably many $p_{2}$. These
are each definable in a $\Sigma_{2}$ way over $L_{\zeta_{2}}$, in the
parameter $(T^{*} ( \varnothing ))_{p_{2}}$. This argument uses that $(T^{*} (
\varnothing ))_{p_{2}} \in L_{\zeta_{1}} \prec_{\Sigma_{1}} L_{s_{1}}$. Let
$T' (p_{2} ) \in L_{\zeta_{2}}$ be {\itshape{II}}'s non-losing quasi-strategy
for $G (A,T^{*} (p_{2} ))$, this time with ``$y=T' (p_{2} )$''$\in
\Delta_{2}^{L_{\zeta_{2}}} (\{T^{*} (p_{2} )\})$. \ (Again these will satisfy
the same definitions as over $L_{\zeta_{m}}$ for any $m \geq 2$.) Note that we
may assume that the countably many trees $T' (p_{2} )$ appear boundedly below
$\zeta_{2}$ (using the $\Sigma_{2}$-admissibility of $\zeta_{2}$). Again for
$p_{3} \in T^{*} (p_{2} )  $ any position of length 3, let $\tau (p_{3} )  $be
some arbitrary but fixed move in $T' (p_{2} )$. Now we consider appropriate
moves $p_{4}$ of length 4, and reapply the lemma with $B=A_{2}$ and $(T^{*}
(p_{2} ))_{p_{4}}$. \ Continuing in this way we obtain a strategy $\tau$ for
{\itshape{II}}, so that $\tau \upharpoonright^{    [1,2k+2)} \omega , $ for
$k< \omega$, is defined by a length $k$-recursion that is
$\Sigma_{2}^{L_{\zeta_{k}}} (\{T\})$.

As the argument continues more and more of the strategy $\tau$ is defined
using successive $( \zeta_{m} ,s_{m} )$ to justify the existence of the
relevant trees in $L_{\zeta_{m}}$. {{\em Knowing\/}} that the trees are there
for the asking, we see that $\tau$ can actually be defined by a
$\Sigma_{2}$-recursion over $L_{\beta_{0}}$ in the parameter $T$ in precisely
the manner given above (the $\Sigma_{2}$-inadmissibility of $\beta_{0}$
notwithstanding).

If $x$ is any play consistent with $\tau , $then for every $n$, by the
defining properties of $T^{*} (p_{2n} )  $given by the relevant application of
the lemma, $x \in \lceil T^{*} (x \upharpoonright 2n) \rceil \subseteq \neg
A_{n.}   $Hence $x \notin A$, and $\tau$ is a winning strategy for{\itshape{
II}} as required. \ Thus $\beta_{0} \geq \gamma_{0}$ is demonstrated.

$( \beta_{0} \leq \gamma_{0} )  $: suppose $\beta_{0} > \gamma_{0}$. Then,
since the existence of a winning strategy for a player in any particular
$\game \Sigma^{0}_{3}$ game would be part of the theory $T^{1}_{\beta_{0}}
=T^{1}_{\alpha_{0}}$ where $\alpha_{0}$ is least with $L_{\alpha_{0}}
\prec_{\Sigma_{1}} L_{\beta_{0}}$, and since moreover that the existence of a
stage $\gamma_{0}$ over which {{\em all\/}} such games have strategies,
amounts also to an existential statement, we have that $\gamma_{0} <
\alpha_{0}$. But this is an immediate contradiction: find a $\psi \in
T^{1}_{\alpha_{0}}$ with $\gamma_{0} < \alpha_{\psi} < \alpha_{0}$. But as
before {\ptwo} has as winning strategy $\sigma$ to play a code for
$L_{\alpha_{\psi}}$. Hence as $\gamma_{0} < \alpha_{\psi}$ such a strategy and
so such a code can be found in $L_{\alpha_{\psi}}$; but again as before, this
contradicts Tarski. Contradiction. Hence $\beta_{0} \leq \gamma_{0}$.
\\ {\hspace*{\fill}} {\qed} Theorem \ref{gamma=beta}

\begin{remark}{\em{
  We make some definitions from the $( \beta_{0} \geq \gamma_{0} )$ part of
  the last proof for later use. We have our starting tree $T$, and the tree of
  non-losing positions for {\ptwo}, $T'$. We shall call these the trees of
  {{\em depth 0\/}}. Then for any $p \in T'$ we argued that $p$ was good, and,
  since $\emp$ was good, we could define the tree $T^{\ast} \left( \emp
  \right)$ - the $L$-least tree witnessing this fact, and thence we had $T'
  \left( \emp \right)$ the tree of non-losing positions for {\ptwo} in \ $G
  (A,T^{*} ( \varnothing ))$. We give the trees $T^{*} ( \emp ) ,T' ( \emp )$
  {{\em depth 1\/}}. Then for any position $p_{1} \in T$ with
  $\ensuremath{\operatorname{lh}} (p_{1} ) =1$, we let $\tau (p_{1} )  $be the
  numerically least move in $T' ( \varnothing )$. \ We call any play, $p_{2}$
  say, of length 2 consistent with this definition of $\tau$ so far, {{\em
  relevant (of length 2). \/}}We wished to apply the lemma again with
  $B=A_{1}$ replacing $B=A_{0}$ and with $(T^{*} ( \varnothing ))_{p_{2}}$
  replacing $T$. We shall call a tree of the form $(T^{*} ( \varnothing
  ))_{p_{2}}$ or $( (T^{*} ( \varnothing ))_{p_{2}} )'$ (the latter the tree
  of non-losing moves for {\ptwo} in $G ( A;  (T^{*} ( \varnothing
  ))_{p_{2}}$)) {{\em relevant trees of depth 1\/}}. \ We then used $(
  \zeta_{2} ,s_{2} )$ to define the $T^{*} (p_{2} )$ (one tree for each
  relevant $ p_{2} )$ and thence the trees $T' (p_{2} )$ to be
  {\itshape{II}}'s non-losing quasi-strategy for $G (A,T^{*} (p_{2} ))$. We
  give trees of the form $T^{*} (p_{2} ) ,T' (p_{2} )$ {{\em depth 2\/}}. For
  $p_{3} \in T^{*} (p_{2} )  $ any position of length 3, $\tau (p_{3} )  $was
  the numerically least move in $T' (p_{2} )$. Again we call such $p_{4}
  =p_{3} \smallfrown \tau (p_{3} )$ {{\em relevant\/}}, and the corresponding
  trees $(T^{*} (p_{2} ))_{p_{4}}$ and $ (T^{*} (p_{2} ))_{p_{4}} )'  $ {{\em
  relevant trees of depth 2\/}}. \ $T^{*} (p_{4} ) ,T' (p_{4} )$ will be of
  {{\em depth 3\/}}. And so forth.}}
\end{remark}

\begin{definition}
  Let $\mathbbm{T}^{k}$ denote the set of trees, and relevant trees, of depth
  $k$, as just defined for $k< \omega$.
\end{definition}

We return now to considering the complexity of $\Game \Sigma^{0}_{3}$.

\begin{theorem}
  \label{GSigma}Let $\alpha_{0}$ be least with $T^{1}_{\alpha_{0}}
  =T^{1}_{\beta_{0}}$ (thus $ \alpha_{0} = \min  S^{1}_{\beta_{0}} )$.
  
  (i) \ $T^{1}_{\alpha_{0}}$ is a complete $\Game \Sigma_{3}^{0}$ set of
  integers.
  
  (ii) Hence the reals of $L_{\alpha_{0}}$ are all $\game \Sigma^{0}_{3}$ sets
  of integers.
\end{theorem}

{\nod}{\pf} The argument is really close to that of the Corollary 2 of
{\cite{We11}}. Indeed there we showed that the $T^{1}_{\alpha_{\psi}}$ (which
occurred cofinally in $L_{\alpha_{0}}$) were $\game \Sigma^{0}_{3}$ sets. Some
details of this are repeated. First remark that (ii) is immediate given (i)
since all the other reals in $L_{\alpha_{0}}$ are all recursive in
$T^{1}_{\alpha_{0}}$ and \ $\game \Sigma^{0}_{3}$, being a Spector class
({{\em v.\/}} {\cite{Mosch3}}), is closed under recursive substitution. We
define a game $G_{\varphi}^{\ast}$ for $\Sigma_{1}$-sentences $\varphi$.

{{\em Rules for II.\/}}

{\noindent}In this game {\ptwo}'s moves in $x$ must be a set of G{\"o}del
numbers for the complete $\Sigma_{1}$-theory of an $\omega$-model of $
\ensuremath{\operatorname{KP}}+ V=L + ( \neg \varphi \wedge
\ensuremath{\operatorname{Det}}( \Sigma^{0}_{3} ) ) .$

Everything else remains the same {{\em mutatis mutandis\/}}: {\pone}'s Rules
remain the same and his task is to find an infinite descending chain through
the ordinals of {\ptwo}'s model. Note that if $\varphi \in
T^{1}_{\alpha_{0}}$, {\pone} now has a winning strategy: for if {\ptwo} obeys
her rules, and lists an $x$ which codes an $\omega$-model $M$ of this theory,
then $M$ is not wellfounded, and has WFP$(M) \cap
\ensuremath{\operatorname{On}} <  \rho ( \varphi )$ where $\rho ( \varphi )$
is defined as the least $\rho$ such that $\varphi \in T^{1}_{\rho +1}$.
However {\pone} playing (just as {\ptwo} did in the main Theorem 4) can find a
descending chain and win. For we have $\ensuremath{\operatorname{WFP}} (M)
\cap \ensuremath{\operatorname{On}} <  \beta_{0}$ and so the argument goes
through, as there are no infinite depth nestings there. \ On the other hand if
$\varphi \notin T^{1}_{\alpha_{0}}$, {\ptwo} may just play a code for the true
wellfounded $L_{\beta^{+}_{0}}$ with $\beta^{+}_{0}$ the least admissible
above $\beta_{0} +1$, and so win. This shows that $T_{\alpha_{0}}^{1}$ is a
$\Game \Sigma_{3}^{0}$ set of integers. \ \ \ \ \

Now suppose $a \in \Game \Sigma_{3}^{0}$. Then we have some $\Sigma_{3}^{0}$
set $A \sset \omega \times  ^{\omega} \omega$ with $n \in a \,\equi \,I$ has a
winning strategy to play into $A_{a} = \{y \in  ^{\omega} \omega   \mid  (a,y)
\in A\}$. Then $a$ is $\Sigma^{L_{\alpha_{0}}}_{1}$ (since all
$\Sigma^{0}_{3}$-games that are a win for {\pone}, have a winning strategy an
element of $L_{\beta_{0}}$, and thence by $\Sigma_{1}$-elementarity, the
$L$-least such is actually an element of $L_{\alpha_{0}}$ - and we merely have
to search through $L_{\alpha_{0}}  $for it) and thus is recursive in
$T_{\alpha_{0}}^{1}$. Hence \ $T_{\alpha_{0}}^{1}$ is a complete $\Game
\Sigma_{3}^{0}$ set of integers. \ \ \ {\hspace*{\fill}} \ \ {\qed} Theorem
\ref{GSigma} and \ref{2.8}(a) (ii)$\lr$(iii).\\

In conclusion: we saw above that $\alpha_{0}$ was the least $\alpha$ with
$T_{\alpha}^{1} =T^{1}_{\beta_{0}}$. Phrased in other terms, by elementary
constructible hierarchy considerations, this is saying that $\alpha_{0}$ is
the minimum of $S^{1}_{\beta_{0}}$. Hence $L_{\alpha_{0}} \prec_{\Sigma_{1}}
L_{\beta_{0}}$ but for no smaller $\delta$ is $L_{\delta} \prec_{\Sigma_{1}}
L_{\beta_{0}}$. Since the statement ``There is a winning strategy for Player
$I$ in $G (A,T)$'' is equivalent in $\ensuremath{\operatorname{KPI}}$ to a
$\Sigma_{1}$-assertion, if true in $L_{\beta_{0}}$ it is true in
$L_{\alpha_{0}}$. In short for those $\Sigma^{0}_{3}$-games that are wins for
$I$ on trees $T \in L_{\alpha_{0}}$, there are strategies for such also within
$L_{\alpha_{0}}$ itself. For those that are wins for Player {{\em II\/}}, when
not found in $L_{\alpha_{0}}$, these may be defined over $L_{\beta_{0}}$. This
somewhat asymmetrical picture reflects the earlier theorems cited above. \ The
theorems of the next section harmonise perfectly with this.\\

Remark: \ (i) Since $\Game \Sigma_{3}^{0}$ is a Spector class, one will have a
$\Game \Sigma_{3}^{0}$-prewellordering of $T_{\alpha_{0}}^{1}$ as a $\Game
\Sigma_{3}^{0}$ set of integers, of maximal length, here $\alpha_{0}$.

{\color{black} We write down one on $T=T^{1}_{\alpha_{0}}$. Abbreviate
$\Gamma = \Game \Sigma_{3}^{0}$ and $\check{\Gamma} = \Game \Pi_{3}^{0}$. We
need to provide relations $\leq_{\Gamma}$ and $\leq_{\check{\Gamma}}$ in
$\Gamma$ and $\check{\Gamma}$ respectively, so that the following hold:
$$T (y) \Longrightarrow \all x \{[T(x) \wedge \rho
(x) \leq \rho (y)]  \Longleftrightarrow x \leq_{\Gamma} y \Longleftrightarrow
x \leq_{\check{\Gamma}} y\}.$$

For the relation $x \leq_{\Gamma} y$, we define the game where {\ptwo}
produces a model $M^{I \!\! I}$ of $T (y) \wedge ( \neg T(x)  \vee \rho (x)
{\nleq} \rho (y))$ and $I$ tries to demonstrate that it is illfounded.
Assume then $T (y)$. If $T (x) \wedge \rho (x) \leq \rho (y)$ then {{\em
either\/}} $( \neg T(x))^{M^{I I}}$ and thus $M^{I I}$ is illfounded with
$\ensuremath{\operatorname{WFP}} (M^{I \!\! I} ) \cap
\ensuremath{\operatorname{On}}< \rho (x)$ and hence $I$ can win as in this
region there are no $\omega$-nested sequences. {{\em Or\/}}: $( \rho (x)
{\nleq} \rho (y))^{M^{I \!\! I}}$. Thus $( \rho (x) > \rho (y))^{M^{I \!\!
I}}$ and again this implies $\ensuremath{\operatorname{WFP}} (M^{I \!\! I} )
\cap \ensuremath{\operatorname{On}}< \rho (x)$ with $I$ winning.

Conversely suppose $x \leq_{\Gamma} y$. Since $T (y)$ is assumed, if $\neg T
(x)$, then {\ptwo} can play a wellfounded model with $(y \wedge \neg x)^{M^{I
I}}$ and win. If $\rho (x)  > \rho (y)$ then again the same can be done. This
proves the first equivalence above. The second is similar, with now $I$
producing a model $M^{I}$ of $T (x) \wedge \rho (x) \leq \rho (y)$ and {\ptwo}
finding descending chains. We leave the details to the reader.}

(ii) One may also write out directly the theories $T^{1}_{\alpha}$ for $\alpha
< \alpha_{0}$ in a $\Game \Pi_{3}^{0}$ form. This should not be surprising: a
$\Game \Sigma_{3}^{0}$ norm as above should have `good' $\Delta ( \Game
\Sigma_{3}^{0} )$ initial segments.

(iii) \ For any set $A \in \Game \Pi_{3}^{0} \backslash \Game \Sigma_{3}$
there will be $n \in A$ so that the winning strategy witnessing this is
definable over $L_{\beta_{0}}$ but not an element thereof. (Otherwise an
admissibility and $\Sigma_{1}$-reflection argument shows that there is a level
$L_{\delta}$ with $\delta < \alpha_{0}$ containing strategies for both $A$ and
its complement. But that would make $A \in \Delta ( \Game \Sigma_{3}^{0} )$ -
a contradiction.)

\begin{corollary}
  \label{sigma=alpha} $\eta_{0} = \alpha_{0}$.
\end{corollary}

{\pf} Since $\Game \Sigma_{3}^{0}$ is a Spector class, and we see that a
complete $\Game \Sigma_{3}^{0}$ set has a $\Game \Sigma_{3}^{0}$ -norm of
length $\alpha_{0}$, standard reasoning shows that there is a $\Game
\Sigma_{3}^{0}$-monotone operator whose closure ordinal is $\alpha_{0}$. Hence
$\eta_{0} = \alpha_{0}$. {\hspace*{\fill}}{\qed}

Results of Martin in {\cite{Martin80}} show that for a co-Spector class,
$\check{\Gamma}$ say, the closure ordinal of monotone
$\check{\Gamma}$-operators, $\omicron (
\check{\Gamma}$-$\ensuremath{\operatorname{mon}}) \dfs \sup \{\omicron ( \Phi
)   \mid \Phi \in \check{\Gamma}$, $\Phi$ \ monotone$\}$, is {{\em
non-projectible\/}}, that is $L_{\omicron \left( \check{\Gamma} \text{-}
\ensuremath{\operatorname{mon}} \right)} \models
\Sigma_{1}$-$\ensuremath{\operatorname{Sep}}.$ \ Moreover \ $o$($\Gamma$) $<$
$\omicron ( \check{\Gamma}$-$\ensuremath{\operatorname{mon}}$ ).

He shows:

\begin{theorem}
  (Theorem D {\cite{Martin80}}) Let $\Gamma$ be a Spector pointclass. (i)
  Suppose that for every $X \sset \omega$, and every $\check{\Gamma} (X)$
  monotone $\Phi$, that $\Phi^{\infty} \in \check{\Gamma} (X)$, then \ $o(
  \check{\Gamma}$-$\ensuremath{\operatorname{mon}}$) is non-projectible, that
  is $S^{1}_{o ( \check{\Gamma} \text{-} \ensuremath{\operatorname{mon}})}$ is
  unbounded in $o( \check{\Gamma}$-$\ensuremath{\operatorname{mon}}$).
  
  (ii) ({{\em \/}}from the proof of his Lemma D.1) \ \ \ \ $o ( \Gamma $-$
  \ensuremath{\operatorname{mon}} )   \in S^{1}_{o ( \check{\Gamma} \text{-}
  \ensuremath{\operatorname{mon}})}$.
\end{theorem}

(He shows too that for Spector classes such as $\game \Sigma^{0}_{3}$, the
supposition of (i) is fulfilled.) If we set $\pi_{0}$ to be the closure
ordinal of $\game \Pi^{0}_{3}$-mon. operators, then in this context we have an
upper bound for $\pi_{0}$:

\begin{lemma}
  $\alpha_{0} < \pi_{0} \leq \beta_{0}$. 
\end{lemma}

{\pf}By (ii) of the last theorem, $\alpha_{0} \in S^{1}_{\pi_{0}}$. But for no
$\beta' > \beta_{0}$ do we have $L_{\alpha_{0}} \prec_{\Sigma_{1}} L_{\beta'}$
(as there are games with winning strategies (for {\ptwo}) in $L_{\beta_{0}
+1}$ for which there are none in $L_{\beta_{0}}$).

{\hspace*{\fill}} \ {\qed}

{{\em Question:\/}} \ Is $\pi_{0} = \beta_{0}$?

\section{Recursion in $\ensuremath{\operatorname{eJ}}$}

\subsection{Kleene Recursion in higher types}

We take some notation and discussion from Hinman {\cite{Hi78}}. \ There was
developed the basic theory of higher type recursion based on an equational
calculus defined by Kleene and refined by him and Gandy in the 1960's. \ The
basic intuition was to define recursions using not just recursive functions on
integers but also allowing recursive schemes using `computable' functions $f:
\omega \times \bai \imp   \omega$ (and similarly for domains which are product
spaces of this type). \ A basic result in this area is that the functions
recursive in $E$ (defined below) are precisely those recursive in $J$, the
`ordinary Turing jump', where we set
$$J ( e,\ensuremath{\boldsymbol{m}}, \ensuremath{\boldsymbol{x}} ) \
 = \left\{ 
 \begin{array}{ll} 0 &\quad \mbox{ if } \{ e \} ( \ensuremath{\boldsymbol{m}} ,
\ensuremath{\boldsymbol{x}} ) \da \\
  1 & \quad \mbox{ otherwise. } 
  \end{array}
  \right.
$$

(We shall follow mostly Hinman in using boldface notation, early or
mid-alphabet roman for integers, but end alphabet roman for elements of \
$\bai$, \ to indicate an (unspecified) number of variables of the given type
in an appropriate product space $^{k} \omega \times  ^{l} \left( \bai \right)$
- which he abbreviates as $^{k,l} \omega$.) Then $E$ (often written $^{2} E$)
is the functional:
$$
E ( x )  = \left\{
 \begin{array}{ll} 0 & \mbox{  if } \ex n ( x ( n ) =0 ) ;\\
 1 & \mbox{ otherwise.}
 \end{array} \right.
 $$

For a fixed type-2 functional $I$ of the kind above - thus a function $I: ^{k}
\omega \times  ^{l} \left( \bai \right) \imp \omega$ such as $E$ or $J$ just
defined, an inductive definition of a set, $\Omega ( I )$, consisting of
equational clauses can be built up in $\omega_{1}$-steps. This defines the
class of those functions $\{ e \}^{I}$ that are recursive in $I$. Of course
such include partial functions, as a descending chain of subcomputation calls
in the tree of computations represents divergence. Just as the clauses of the
induction and the set $\Omega ( I )$ is an expansion of those clauses and
functions of type-1 recursion, also due to Kleene and yielding an inductive
set $\Omega$, we shall wish to expand the notion of `computation' further
along another axis.

Our notation for computation will be modelled on that of the transfinite
machine model, the `infinite time Turing machine' introduced by Hamkins and
Kidder {\cite{HL}}. The signifying feature of such ITTM's is the transfinite
number of stages that they are allowed to run their standard finite Turing
program, on their one-way infinite tape. The behaviour at limit stages is
defined by a `liminf' rule for the cell values of $0$ or $1$, and a replacing
of the read/write head back at the start of the tape, and finally a special
`limit state' $q_{L}$ is entered into.

Actually the formalism is quite robust: one may change details of these
arrangements without altering the computational power. In {\cite{HL}} they
considered a 3-tape arrangement (for Input, Scratch Work, and Output). The
paper {\cite{HaSe}} shows this can be reduced to 1-tape (if the alphabet has
more than two symbols!). One can change the limit behaviour so that instead of
a liminf value being declared for each cell's value, it simply becomes blank -
for ambiguity - if it has changed value cofinally in the limit stage (Theorem
1 of {\cite{W8}}). Similarly the special state $q_{L}$ is unnecessary: one may
define the ``next instruction'' at a limit stage to be the instruction, or
transition table entry, whose number is the liminf of the previous instruction
numbers - this has the machine entering the outermost subroutine that was
called cofinally in the stage. Likewise the Read/Write head may be placed at
the cell numbered according to the liminf values of the cells visited prior to
that limit stage (unless that liminf is now infinite, in which case we do
return the head to the starting cell). \ All of these variants make no
difference to the functions computed.

We shall review the following facts related to such machines.

\subsection{Infinite Time Turing Machine computation}

Such ITTM's have two modes of producing results: a program can halt outright
with an infinite string of $0,1$'s on the part of the tape designated for
output (the `output tape') but it may also have some `eventual output': the
contents of the output tape may have settled down to a fixed value, whilst the
machine is still churning away perhaps moving its head around and fiddling
with the scratch portion of the tape. Nevertheless on a given fixed input
(some $x \in \can$ may be written to a designated portion of the tape, the
`input tape') any ITTM machine will eventually start to cycle - and by the
starting point of that cycling, designated $\zeta ( x )$, if the output
settles down, then it will have done so by $\zeta ( x )$.

This last feature is in fact, quite fundamental for the study of ITTMs. We may
regard a machine $P_{e} ( x )$ in this context, as having come to a conclusion
- the contents of the output tape - but has not formally reached a halting
state in the usual sense.

\begin{definition}
  We shall say that a computation $P_{e} ( x )$ is {{\em convergent\/}} {{\em
  to\/}} $y $(and write $P_{e} ( x ) | y )$ if it enters a halting state in
  the usual sense, or if it has eventually settled output. We shall say that
  ``$y$ {{\em is (eventually)-ittm-recursive in \/}}$x$''. If it does not have
  settled output, we shall write $P_{e} ( x ) \ua$.
\end{definition}

This enshrines our taking {{\em (eventually) settled output\/}}, as the
criterion of a successful computation. We shall be interested in eventual
output of this sort, as well as the more restricted strictly halting variety.
Both types of computation, the usual halting, and the `eventually constant'
output tape outlined above, we shall regard, and term, as `convergent' -
thinking of `halting' as only a special kind of eventually settled output.
Given a set $A \sset \omega \cup \can$, this can be used as an oracle for an
ITTM in a familiar way: ?{{\em Is the integer on (or is the whole of) the
current output tape contents an element of $A$?\/}} and receive a $1/0$ answer
for ``Yes''/``No''. We identify elements of $\omega$ as coded up in {\can} in
some fixed way, and so may consider such $A$ as subsets of $\can$. But
further: since having $A$ respond with one $0/1$ at a time can be repeated, we
could equally as well allow $A$ to return an element $f \in \can$ as a
response (we have no shortage of time). We could then also allow as
functionals also $A: \can \imp \can$. However for the moment we shall only
consider functionals into $\omega$. Some examples follow.

\begin{definition}
  (The infinite time jump $\ensuremath{\operatorname{iJ}}$)\\
  (i)  We write 
 $  \{ e \} (
\ensuremath{\boldsymbol{m}} , \ensuremath{\boldsymbol{x}} ) \da  $
  if  the
$e$'th ittm-computable function with  input   $\ensuremath{\boldsymbol{m}},
 \ensuremath{\boldsymbol{x}}  $ has a halting value.\\
 (ii) We then define $iJ$ by: 
$$\ensuremath{\operatorname{iJ}} ( e,\ensuremath{\boldsymbol{m}},
\ensuremath{\boldsymbol{x}} )\ \ =  
\left\{
\begin{array}{ll}
 1 & \mbox{ if } \{ e \} (
\ensuremath{\boldsymbol{m}} , \ensuremath{\boldsymbol{x}} ) \da  ; \\
 0 & \mbox{ otherwise.}
 \end{array}
\right.
$$
\end{definition}

\begin{definition}
  (The eventual jump $\ensuremath{\operatorname{eJ}}$) \\
  (i)  We write 
 $  \{ e \} (
\ensuremath{\boldsymbol{m}} , \ensuremath{\boldsymbol{x}} ) |  $
  if  the
$e$'th ittm-computable function with  input   $\ensuremath{\boldsymbol{m}},
 \ensuremath{\boldsymbol{x}}  $ has an eventually settled value.\\
 (ii) We then define $eJ$ by: 
$$\ensuremath{\operatorname{eJ}} ( e,\ensuremath{\boldsymbol{m}},
\ensuremath{\boldsymbol{x}} ) \ \ = 
\left\{ 
\begin{array}{ll} 1 & \mbox{ if }\{ e \} (
\ensuremath{\boldsymbol{m}} , \ensuremath{\boldsymbol{x}} ) | ;\\
 0 & \mbox{ otherwise 
(for which we write $\{ e \} ( \ensuremath{\boldsymbol{m}} ,
\ensuremath{\boldsymbol{x}} ) \ua$).}
\end{array}
\right.
$$
\end{definition}

{\nod}These are both total functionals. We shall be interested in functions
recursive in $\ensuremath{\operatorname{eJ}}$. But first we summarise some
facts about ordinary ittm's.\\

\nod{\textbf{Fact 1}} {\cite{HL}} shows:

(i) That $\Pi^{1}_{1}$-predicates are decidable: given a code $x \in \cant$,
there's an ittm that will decide whether $x \in
\ensuremath{\operatorname{WO}}$ or not.

(ii) There's a program number $e$ so that $P_{e} ( x )$ will halt with a code
for $( L_{\alpha} , \in )$ if $x \in \ensuremath{\operatorname{WO}} \wedge || x
|| = \alpha  $.

(iii) For $z \in \cant$, the set of {{\em ittm-writable-in-z\/}} reals, is
the set $\mathcal{W}^{z} \sset \cant$ where

$\mathcal{W}^{z} = \left\{ x \in \cant \mid \ex e P_{e} ( z ) \right.$ halts
with output $ x \}$.

(iv) \ The set of {{\em ittm-eventually-writable-in-$z$\/}} reals, is the set

$\mathcal{E}\mathcal{W}^{z} = \left\{ x \in \cant \mid \ex e  ( P_{e} ( z ) 
\right.$ {{\em has $x$ written on its output tape from some point in time
onwards\/}}$  ) \}$.\\

\nod{\textbf{Fact 2}} {\cite{W09}} shows:

(i) Let $( \zeta , \Sigma )$ be the lexicographically least pair of ordinals
so that $L_{\zeta} \prec_{\Sigma_{2}} L_{\Sigma}$. Let $\lambda$ be the least
ordinal with $L_{\lambda} \prec_{\Sigma_{1}} L_{\zeta}$. \ \ Then ({{\em The
``$\lambda$-$\zeta$-$\Sigma$-Theorem''\/}}), $L_{\lambda} \cap \cant
=\mathcal{W}$, $L_{\zeta} \cap \cant =\mathcal{E}\mathcal{W}$. As is easily
seen all three ordinals are limits of $\Sigma_{2}$-admissibles, whilst
$\lambda$ is $\Sigma_{1}$- but not $\Sigma_{2}$-admissible, and $\Sigma$ is
not admissible at all.

(ii) (a) Any computation $P_{e} ( n )$ that halts (in the usual sense) does so
by a time $\alpha < \lambda$.

(b) Any computation $P_{e} ( n )$ that eventually has a settled output tape,
does so by a time $\alpha < \zeta$. \

(c) Both $\lambda$ and $\zeta$ are the suprema of such fully ``halting''
times, and ``eventual convergence'' times, over varying $e,n \in \omega$,
respectively.

(iii) $T^{1}_{\lambda} \equiv_{1} h$, and $T^{2}_{\zeta} \equiv_{1} \tilde{h}$
where $h= \{ e \mid P_{e} ( e ) $ {{\em reaches a halting state\/}}$ \}$ and
$\tilde{h} = \{ e | P_{e} ( e )  $ {{\em eventually has settled output\/}}$
\}$.

(iv) It is a consequent of $( \ensuremath{\operatorname{iii}} )$ that a {{\em
universal machine\/}} (on integer input) has {{\em snapshots\/}} of its
behaviour which, when first entering a final loop at stage $\zeta$, will
repeat with the same snapshot at time $\Sigma$; moreover (1-1) in those
snapshots is the theory $T^{2}_{\zeta}$.

(v) \ {{\em Recursion,\/}} and $\ensuremath{\operatorname{Snm}}$ {{\em
Theorems\/}} may be proved in the standard manner ({\cite{HL}}); there are
appropriate versions of the Kleene {{\em Normal Form Theorems\/}} \
({\cite{W09}}).\\

The usual argument shows:

\begin{theorem}
  {{\em (The $\tmop{eJ}$-Recursion theorem)\/}} \ If $F ( e,
  \ensuremath{\boldsymbol{m,x}} )  $is recursive in
  $\ensuremath{\operatorname{eJ}}$, there is $e_{0} \in \omega$ so that
  
  \ \ \ \ \ \ \ \ \ \ \ \ \ \ \ \ \ \ \ \ \ \ \ \ \ \ \ \ \ \ \ \ \ \ \ \ \ \
  \ \ \ $\varphi^{\ensuremath{\operatorname{eJ}}}_{e_{0}} (
  \ensuremath{\boldsymbol{m,x}} )  =F ( e_{0} , \ensuremath{\boldsymbol{m,x}}
  )$.
\end{theorem}

\subsubsection{More on Extendability}

As Fact 2 (i) above shows, the relation of ``$L_{\zeta}$ has a
$\Sigma_{2}$-extension to $L_{\Sigma}$'' is fundamental to this notion.\\

\nod{\textbf{Fact 2}} (contd.)

(vi) There is moreover a {{\em theory machine\/}} that writes codes for
$L_{\alpha}$ and their $\Sigma_{\omega}$-theories, and hence their
$\Sigma_{2}$-theories, $T^{2}_{\alpha}$, in an {{\em ittm\/}}-computable
fashion for any $\alpha < \Sigma$, uniformly in $\alpha$. If for
$\ensuremath{\operatorname{Lim}} ( \lambda )$ we write $\hat{T}_{\lambda} \dfs
\ensuremath{\operatorname{Liminf}}_{\alpha \rightarrow \lambda  }
T^{2}_{\alpha}$, then there is a uniform index $e \in \omega$ that shows that
$W^{\hat{T}_{\lambda}}_{e} =T^{2}_{\lambda}$, {\ie} $T^{2}_{\lambda}$ is r.e.
in $\hat{T}_{\lambda}$ uniformly in $\lambda$. (See Lemma 2.5 of
{\cite{W2014}}. \ Moreover for those $\lambda$ with $L_{\lambda} \models
\Sigma_{1}$-$\ensuremath{\operatorname{Sep}}$, $T^{2}_{\lambda} =
\hat{T}_{\lambda}$.)

(vii) For the lexicographically least extendible pair $( \zeta , \Sigma )$,
whilst $\omega^{T^{2}_{\zeta}}_{1\ensuremath{\operatorname{ck}}} < \Sigma$, it
is the case that $\lambda ( T^{2}_{\zeta} ) > \Sigma$.\\

We make some further definitions concerning extendability.

\begin{definition}
  (The $\Sigma_{2}$-extendibility tree) We let $( \mathcal{T}, \prec )$ be the
  natural tree on such pairs under inclusion: as follows: \ if $( \zeta' ,
  \Sigma' ) $,$ \overline{( \zeta } , \bar{\Sigma} )$ are any two countable
  $\Sigma_{2}$-extendable pairs, then set \ $( \zeta' , \Sigma' )  \prec 
  \overline{( \zeta } , \bar{\Sigma} )$ iff $\zeta' \leq \bar{\zeta} <
  \bar{\Sigma} < \Sigma'$.
\end{definition}

{\bu} If we had allowed the inequality $\bar{\Sigma} \leq \Sigma'$ rather than
a strict inequalitiy in the last definition we could have defined a larger
relation $\prec'$, and a larger tree $( \mathcal{T}' , \prec' )$; however this
would not have been wellfounded: if $L_{\Sigma} \models
\Sigma_{2}$-$\ensuremath{\operatorname{Sep}}$ then it is easy to see that
$\left( \mathcal{T}' \rest \Sigma +1, \prec' \right)$ is illfounded.

\begin{lemma}
  Let $\delta$ be least such that $L_{\delta} \models \Sigma_{2}$-Sep. ; let
  $\alpha  $be maximal so that $\left( \mathcal{T}' \rest \alpha , \prec'
  \right)$ is wellfounded (where $\ensuremath{\operatorname{Field}} \left(
  \mathcal{T}' \rest \alpha \right) \dfs \{ ( \zeta , \Sigma ) 
  \ensuremath{\operatorname{extendible}} \mid \Sigma < \alpha \}$). Then
  $\delta = \alpha$.
\end{lemma}

{\pf} $( \leq )$: Suppose $\delta > \alpha$. Then \ $\left( \mathcal{T}' \rest
\delta , \prec \right)$ is illfounded. So there is an infinite sequence of
extendible pairs $( \zeta_{n} , \Sigma_{n} )$ with $( \zeta_{n+1} ,
\Sigma_{n+1} ) \subset ( \zeta_{n} , \Sigma_{n} )$. By wellfoundedness of the
ordinals there is an infinite subsequence $( \zeta_{n_{i}} , \Sigma_{n_{i}} )$
with all $\Sigma_{n_{i}}$ equal to a fixed $\Sigma$, whilst $\zeta_{n_{i}} <
\zeta_{n_{i+1}}$. Let $\zeta^{\ast} = \sup_{i}   \zeta_{n_{i}}$. Then we have
$L_{\zeta_{n_{i}}} \prec_{\Sigma_{2}} L_{\zeta_{n_{i+1}}} \prec_{\Sigma_{2}}
L_{\zeta^{\ast}}$. Then $\zeta^{\ast}$ is not $\Sigma_{2}$-projectible, and
hence $L_{\zeta^{\ast}} \models \Sigma_{2}$-Sep. But $\zeta^{\ast} < \delta$.
Contradiction.

($ \geq ) : L_{\delta} \models \Sigma_{2}$-Sep. \ Then $S^{2}_{\delta}$ is
unbounded in $\delta$. Let $\delta_{i} < \delta_{i+1}$ be a cofinal sequence,
for $i< \omega$. Then check that $\la ( \delta_{i} , \delta ) \mid i< \omega
\ra $ is a $\prec$-descending sequence in $\mathcal{T}' \rest \delta
+1$. So $\alpha \leq \delta$. {\hspace*{\fill}}{\qed}\\

For $E $ a class of ordinals, let $E^{\ast}$ denote the class of its limit
points.
\pagebreak
\begin{definition}
  Define by recursion on $0< \alpha \in \ensuremath{\operatorname{On}}$ the
  class $E^{\alpha}$ the class of $\alpha ( $-$ \Sigma_{2} )$-extendible
  ordinals:
  $$  \begin{array}{lcl}
   E^{1} & = &\{ {\zx{1}}\mid   \zx{1} \mbox{ {{\em is extendible but not a limit of
  extendibles\/}}} \}; \\
  E^{\alpha +1} & = & \{ {\zx{{\alpha}}}\mid   \zx{\alpha} \in ( E^{\alpha}
  )^{\ast} \cap E^{0}  \};\\
  E^{\lambda} & = & \bigcap_{\alpha < \lambda} E^{\alpha} \cap E^{0}.\\
  E^{\geq \alpha} & = &\bigcup_{\beta \geq \alpha} E^{\beta} \mbox{ {\etc} }
  \end{array} $$
\end{definition}

Here we decorate the variable $\zeta$ with the prefix indicating its level of
extendability. We shall let $\sx{\alpha}$ indicate that for some
$\zx{\alpha}$, ({\zx{{\alpha}}},{\sx{{\alpha}}}) is an $\alpha$-extendible
pair. \ Note that for any $\gamma$ the least element of $E^{\geq \alpha}$
greater than $\gamma$ is always an element of $E^{\alpha}$, {\ie} is
$\alpha$-extendible.

\subsection{The Lengths of computations}

We analyse the tree of subcomputations to define the notion of {{\em absolute
length\/}} of the linearised {{\em absolute computation\/}} corresponding to
some $P^{I}_{e} ( \ensuremath{\boldsymbol{m}}  ,
\ensuremath{\boldsymbol{\nobracket x )}} $.

\begin{definition}
  The {{\em local length\/}} of a computation $P^{I}_{e} (
  \ensuremath{\boldsymbol{m}}  , \ensuremath{\boldsymbol{\nobracket x )}} $ in
  a type-2 oracle $I$, is the least $\sigma_{0}$ (when defined) so that the
  snapshot at $\sigma_{0}$ is the repeat of some earlier snapshot $\zeta_{0} <
  \sigma_{0}$, and so that the snapshot at $\sigma_{0}$ recurs unboundedly in
  $\ensuremath{\operatorname{On}}$.
\end{definition}

The local length has all the relevant information then in the calculation:
everything thereafter is mere repetition ($\sigma_{0}$ will be undefined if
$P^{I}_{e} ( \ensuremath{\boldsymbol{m}}  , \ensuremath{\boldsymbol{\nobracket
x )}} $ is divergent, that is, has an ill-founded computation tree). Another
description of it is as the ``top level'' length of the computation, which
disregards the lengths of the subcomputation calls below it. We now describe a
computation recursive in the type-2 functional
$\ensuremath{\operatorname{eJ}}$. In fact we give a representation in terms of
ITTM's. $\p{e} ( \ensuremath{\boldsymbol{m}}  ,
\ensuremath{\boldsymbol{\nobracket x )}} $ will represent the $e$'th program
in the usual format with appeal to oracle calls possible. We are thus
considering computation of a partial function
$\Phi^{\ensuremath{\operatorname{eJ}}}_{e} :\mbox{}  ^{k} \omega \times  ^{l}
(^{\omega} 2 ) \rightarrow  ^{\omega} 2$. Such a computation may
conventionally halt, or may go on for ever through the ordinals. The
computation of $\p{e} ( \ensuremath{\boldsymbol{m}}  ,
\ensuremath{\boldsymbol{\nobracket x )}} $ proceeds in the usual ittm-fashion,
working as TM at successor ordinals and taking $\liminf$'s of cell values
{\etc} at limit ordinals. At time $\alpha$ an oracle query may be initiated.
We shall conventionally fix that the real being queried is that infinite
string on the even numbered cells of the scratch type. If this string is $(
f,y_{0} ,y_{1}  \ldots    , )$ then the query is ?{{\em Does $\p{f} ( y )$
have eventually settled output tape\/}}?, and at stage $\alpha +1$ receives a
$1/0$ value corresponding to ``Yes/No'' respectively. We thus regard
$\ensuremath{\operatorname{eJ}}$ as the ``eventual jump'' and intend the
following:
$$\ensuremath{\operatorname{eJ}}= \{ \pa{\langle f,y \rangle ,i} \mid i=1
\mbox{ and } \p{f} ( y ) | \mbox{ or \ } i=0 \mbox{ and } \p{f} ( y ) \ua
\}$$
Here, $\p{f} ( y ) \ua  $ denotes that the computation $\p{f} ( y )$ loops but
has no settled output, it is {{\em not\/}} the notation for a computation
whose tree has an ill-founded branch. (Compare with above for the type-2
recursion in $J$: {{\em divergence\/}} occurs if there is an
illfounded-founded branch in the tree of evaluations.) As is intended, $\p{f}
( y )$ has the opportunity to make similar oracle calls, and we shall thus
have a {{\em tree\/}} representation of calls made. \ We wish to represent the
overall order of how such calls are made, and indeed the ordinal times of the
various parts of the computation as it proceeds. Overall we have a `depth
first' mode of evaluation of a tree of subcomputations. We therefore make the
following conventions. During the calculation of $\p{e} (
\ensuremath{\boldsymbol{m}} , \ensuremath{\boldsymbol{\nobracket x )}} $ (the
topmost node $\nu_{0}$ at Level $0$, in our tree
$\text{$\mathfrak{T}=\mathfrak{T} ( \ensuremath{\boldsymbol{e,m}}  ,
\ensuremath{\boldsymbol{\nobracket x )}} $)}$ let us suppose the first oracle
query concerning $\p{f_{0}} ( y_{0} )$ is made at stage $\delta_{0}$. We write
a node $\nu_{1}$ below $\nu_{0}$, and explicitly allow the computation
$\p{f_{0}} ( y_{0} )$ to be performed at this Level $1$. The `local time' for
this computation, of course starts at $t=0$ - although each stage is also
thought of as one more step in the overall computation of the computation
immediately above: namely $\p{e} ( \ensuremath{\boldsymbol{m}} ,
\ensuremath{\boldsymbol{\nobracket x )}} $. \ Suppose $\p{f_{0}} ( y_{0} )$
makes no further oracle calls and the local length of $\p{f_{0}} ( y_{0} )$ is
$\sigma_{1}$. \ Control, and the correct $1/0$ bit is then passed back up to
Level 0, and the master computation proceeds.

We deem that $\delta_{0} + \sigma_{1}$ steps have occurred so far towards the
final {{\em absolute length\/}} of the calculation $H=H ( e,
\ensuremath{\boldsymbol{m}} , \ensuremath{\boldsymbol{\nobracket x )}} $, of
$\p{e} ( \ensuremath{\boldsymbol{m}} , \ensuremath{\boldsymbol{\nobracket x
)}} $.

However if $\p{f_{0}} ( y_{0} )$ has made an oracle query, let us suppose the
first such was $? \p{f_{1}} ( y_{1} ) ?$, then a new node $\nu_{2}$ is placed
below $\nu_{1}$. If this piece of computation at $\nu_{2}  $takes $\sigma_{2}$
steps without oracle calls, to cycle before control and the result is passed
back up to $\nu_{1}$, ({\ie} the local length of $\p{f_{1}} ( y_{1} )$ is
$\sigma_{2}$) then those $\sigma_{2}$ steps will have to be be part of the
overall length of calculation for $\p{e} ( \ensuremath{\boldsymbol{m}} ,
\ensuremath{\boldsymbol{\nobracket x )}} $ - \ although those $\sigma_{2}$
steps only counted for 1 step in the local length of $\p{f_{0}} ( y_{0} )$'s
calculation. If the $\p{e} ( \ensuremath{\boldsymbol{m}} ,
\ensuremath{\boldsymbol{\nobracket x )}} $ converges then we shall have its
computation tree $\mathfrak{T}=\mathfrak{T} ( e, \ensuremath{\boldsymbol{m}} ,
\ensuremath{\boldsymbol{\nobracket x )}} $, a finite path tree (with
potentially infinite branching) and some countable rank. $\mathfrak{T}$ will
be labelled with nodes $\{ \nu_{\iota} \}_{\iota < \eta ( \mathfrak{T} )}$
that are visited by the computation in increasing order (with backtracking up
the tree of the kind indicated). Thus $\nu_{\iota}$ is first visited only
after all $\nu_{\tau}$ have been visited for $\tau < \iota$. The $\beta$'th
oracle call to Level $k$ will generate a node placed to the right of those so
far at Level $k$ (and thus to the right of those with lesser indices $\alpha <
\beta$ at that level). The tree will thus have a linear leftmost branch,
before any branching occurs.

Just as the Kleene equational calculus can be seen to build up in an
inductive fashion a set of indices $\Omega [ I ]$ for successful computations
recursive in $I$, (see Hinman {\cite{Hi78}}, pp. 259-261) so we can define the
graph of $\ensuremath{\operatorname{eJ}}$ as the fixed point of a monotone
operator $\Delta$ on $\omega \times \omega^{< \omega} \times ( \omega^{\omega}
)^{< \omega} \times 2$.\\

We set $\Delta ( X ) = $: 
$$\begin{array}{l}
\left\{ \langle \pa{e, \ensuremath{\boldsymbol{m}} ,
\ensuremath{\boldsymbol{x}}} ,i \rangle   |  \right. P_{e}^{X} (
\ensuremath{\boldsymbol{m}} , \ensuremath{\boldsymbol{\nobracket x )}} \mbox{{{ \em
is an ittm-computation making only oracle calls}} 
} \\ 
 \,\,  \langle \pa{e' ,
\ensuremath{\boldsymbol{m'}} , \ensuremath{\boldsymbol{x'}}} ,i' \rangle \in
X \mbox{{{ \em with $i=1/0$ if the resulting output is eventually settled or
not$\nobracket \}$.\/}}} 
\end{array}
$$

Let $\Delta^{0} = \emp$; $\Delta^{\alpha +1} = \Delta ( \Delta^{\alpha} )$;
$\Delta^{< \lambda} = \bigcup_{\alpha < \lambda} \Delta^{\alpha}$ \&
$\Delta^{\lambda} = \Delta ( \Delta^{< \lambda} )$ in the usual way. Then the
least fixed point of $\Delta$ is the function
$\ensuremath{\operatorname{eJ}}$.

\begin{definition}
  With $\ensuremath{\operatorname{eJ}}$ as just defined:
  
  $\p{e} ( \ensuremath{\boldsymbol{m}} , \ensuremath{\boldsymbol{\nobracket x
  )}} $ is {{\em convergent\/}} if $\pa{e, \ensuremath{\boldsymbol{m}} ,
  \ensuremath{\boldsymbol{x}}} \in \ensuremath{\operatorname{dom}} (
  \ensuremath{\operatorname{eJ}} )$. \ Otherwise it is {{\em divergent.\/}}
\end{definition}

Assuming $\p{e} ( \ensuremath{\boldsymbol{m}} ,
\ensuremath{\boldsymbol{\nobracket x )}} $ convergent, we may define by
recursion a function $H ( f_{i} ,y_{i} )$ for $1 \leq \iota < \eta (
\mathfrak{T} )$, giving that absolute length of the calculation at node
$\nu_{\iota}$ taking into account the computations at nodes below it. Suppose
the oracle queries made by $\p{e} ( \ensuremath{\boldsymbol{m}} ,
\ensuremath{\boldsymbol{\nobracket x )}} $ \ at Level 0, were
$\p{f_{\iota_{j}}} ( y_{\iota_{j}} )$ for $j< \theta$, and they were made at
increasing local times $\delta_{j}$ for $j< \theta$ in $\p{e} (
\ensuremath{\boldsymbol{m}} , \ensuremath{\boldsymbol{\nobracket x )}} $, then
let $\overline{\delta_{i}}$ be defined by:
$$
\begin{array}{lcl}
\overline{\delta_{0}} & = & \delta_{0} ;\\
\overline{\delta_{}}_{j+1} & = &\delta_{j+1} - \delta_{j};\\
\bar{\delta}_{\lambda} & = &\delta_{\lambda} - \sup \{ \delta_{k}   \mid  k< \lambda \} 
 \end{array}
$$
{\nod}then the {{\em absolute length\/}} of the calculation is the wellordered
ordinal sum:

$$
H ( e, \ensuremath{\boldsymbol{m}} , \ensuremath{\boldsymbol{\nobracket
x )}} = 
\left\{\begin{array}{ll}
\sum^{\theta}_{j=0} ( \bar{\delta}_{j} +H ( f_{\iota_{j}}
,y_{\iota_{j}} ) ) & \mbox{ if } \theta >0;\\
\Sigma (
\ensuremath{\boldsymbol{x}} ) &\mbox{ otherwise}.
\end{array} \right.
$$

{\nod}of course assuming by induction that the absolute lengths of the
computations $H ( f_{\iota_{j}} ,y_{\iota_{j}} )$ have been similarly defined.

We call the master computation $\p{e} ( \ensuremath{\boldsymbol{m}} ,
\ensuremath{\boldsymbol{\nobracket x )}} $ together with all the
subcomputations of the tree explicitly performed, the {{\em absolute
computation\/}} (as opposed to the top level `local computation' with simple
$1$-step queries). \

{\bu} It is possible, and easy, to design an index $f \in \omega$, \ so that
$\p{f} ( 0 )$ has absolute length $H \left( f,0, \emp \right)$ greater than
the looping length of the top level computation. \ Hence for performing a
computation together with all its subcomputations as a tree, and seeing how
the absolute computation relates to extendability in the $L$ hierarchy, this
has to be done in suitably large admissible sets.

\begin{lemma}
  Suppose $\p{e} ( \ensuremath{\boldsymbol{m}} ,
  \ensuremath{\boldsymbol{\nobracket x )}} $ is a convergent computation with
  tree $\mathfrak{T} \in M$, and with $\ensuremath{\boldsymbol{x}} \in
  M,$where $M$ is a transitive admissible set. Let $\theta
  =\ensuremath{\operatorname{On}}^{M}$. Suppose for every node $\nu_{\iota}$
  in $\mathfrak{T}$ that the computation at the node $\p{f_{\iota}} (
  y_{\iota} )$ has {{\em local \/}}length $\psi_{\iota} < \theta$ (this
  includes the local length of $\p{e} ( \ensuremath{\boldsymbol{m}} ,
  \ensuremath{\boldsymbol{\nobracket x )}} $, being at $\nu_{0}$, is some
  $\psi_{0} < \theta$). Then $H ( e, \ensuremath{\boldsymbol{m}} ,
  \ensuremath{\boldsymbol{\nobracket x )}} < \theta $.
\end{lemma}

{\pf} The required ordinal sum can be performed by an induction on the {{\em
rank\/}} of the nodes in the tree, setting $0=\ensuremath{\operatorname{rank}}
( \nu_{\iota} )$, for those $\iota$ with $\nu_{\iota}$ a terminal point of a
path leading downwards from $\nu_{0}$. This can be effected inside the
admissible set $M$. {\hspace*{\fill}}{\qed}

Better:

\begin{lemma}
  Suppose $\p{e} ( \ensuremath{\boldsymbol{m}} ,
  \ensuremath{\boldsymbol{\nobracket x )}} $ is a convergent computation with
  its computation tree $\mathfrak{T} \in M$, and with
  $\ensuremath{\boldsymbol{x}} \in M,$where $M$ is a transitive admissible
  set, closed under the function $x \rightarrowtail \tilde{x}$. Let $\theta
  =\ensuremath{\operatorname{On}}^{M}$. Then $H ( e,
  \ensuremath{\boldsymbol{m}} , \ensuremath{\boldsymbol{\nobracket x )}} <
  \theta $.
\end{lemma}

{\pf}This is similar to the above. By induction on
$\ensuremath{\operatorname{rk}} \left( \gott \right) = \eta < \theta$. Note
first that the closure of $M$ ensures that for all $y \in M$, that $\Sigma ( y
) < \theta$. Suppose true for all such trees of convergent computations $P_{f}
( y )$ of smaller rank than $\eta$, for $y \in M.$ Suppose $\p{e} ( x )$ makes
queries at {{\em local\/}} times $\pa{\delta_{i} \mid i< \tau}$ to nodes at
Level 1. Note that $\tau < \theta$ as {\gott}$\in M$. Suppose the calls are to
the subtrees $\pa{\gott_{i} \mid i< \tau}$ with $( f_{i} ,y_{i} )$ passed down
at time $\delta_{i}$ and $\overline{y_{i}}$ is the real passed up at local
time $\delta_{i} +1$. Let the snapshot at Level 0 at time $\gamma$ be $s (
\gamma )$. (Thus we assume $s ( \delta_{i} +1 )$ contains the information of
$\overline{y_{i}}$.) Now notice that $\delta_{0} < \Sigma ( x )$ (because the
computation prior to $\delta_{0}$ is (equivalent to) an ordinary ittm
computation, which of course eventually converges at time $\Sigma ( x )$.) If
we set
$$
\begin{array}{lclrcl}
\overline{\delta_{0}} & = &\delta_{0} ; & s_{0} & = & x\\
\overline{\delta_{}}_{j+1} & =  &\delta_{j+1} - \delta_{j}; & s_{j+1} & = &s (
\delta_{j} )\\
\bar{\delta}_{\lambda} & = & \delta_{\lambda} - \sup \{ \delta_{k}   \mid  k<
\lambda \}  ; & s_{\lambda} & = &s ( \sup \{ \delta_{k}   \mid  k< \lambda \} );
\end{array}
$$

{\nod}then $\bar{\delta}_{j+1} < \Sigma ( s_{j} )$ (as the time to the next
query, if it exists, is always less than the least $s_{j}$-2-extendible by the
same reasoning). Similarly $\bar{\delta}_{\lambda} < \Sigma ( s_{\lambda} )$.
By assumption on $M$, all such $\Sigma ( s_{j} )$ are less than $\theta$. \
Consequently if $H ( f_{i} ,y_{i} ) = \theta_{i} < \theta$, the whole length
of the computation is bounded:

$$H ( e,
\ensuremath{\boldsymbol{m}} , \ensuremath{\boldsymbol{\nobracket x )}} =
\sum_{i=0}^{< \tau} \overline{\delta_{}}_{i} + \theta_{i}   \leq  
\sum_{i=0}^{< \tau} \Sigma ( s_{i} ) + \theta_{i}  <  \theta .$$

{\hspace*{\fill}}{\qed}

\begin{definition}
  (i) The {{\em Level\/}} of the computation $\p{e} (
  \ensuremath{\boldsymbol{m}} , \ensuremath{\boldsymbol{\nobracket x )}} $ at
  time $\alpha <H ( e, \ensuremath{\boldsymbol{m}} ,
  \ensuremath{\boldsymbol{\nobracket x )}} $, denoted $\Lambda ( e, (
  \ensuremath{\boldsymbol{m}} , \ensuremath{\boldsymbol{\nobracket x )}} ,
  \alpha ) $, is the level of the node $\nu_{\iota}$ at which control is based
  at time $\alpha$, where:
  
  (ii) \ the {{\em level\/}} of a node $\nu_{\iota}$ is the length of the path
  in the tree from $\nu_{0}$ to $\nu_{\iota}$.
\end{definition}

Thus for a convergent computation, at any time the level is a finite number
(`depth' would have been an equally good choice of word). A divergent
computation is one in which $\mathfrak{T} ( e, \ensuremath{\boldsymbol{m}} ,
\ensuremath{\boldsymbol{\nobracket x )}} $ becomes illfounded (with a
rightmost path of order type then $\omega$).

\begin{lemma}
  The computation $\p{e} ( x )$ converges if and only if there exists some
  $x$-$\Sigma_{2}$-extendible pair $( \zeta , \Sigma )$ so that $\Lambda (
  e,x, \zeta ) =0$.
\end{lemma}

{\pf} Suppose $\p{e} ( x ) | $. If $\pe ( x ) \da$ conventionally then the
conclusion is trivial as then for all sufficiently large
$x$-$\Sigma_{2}$-extendible pairs $( \zeta , \Sigma )$, the machine has halted
at Level 0. If otherwise, then the computation $\p{e} ( x )$ will loop forever
through the ordinals. But, using the definition of the $\liminf$behaviour at
limit stages, it is easy to argue that there is a cub subset $C \sset
\omega_{1}$ of points $\alpha , \beta$ with the snapshots of the computation
at these times identical, and with $\Lambda ( e,x, \alpha ) = \Lambda ( e,x,
\beta ) =0$. Now find a pair $( \zeta , \Sigma' )$ both in $C$, with
$L_{\zeta} [ x ] \prec_{\Sigma_{2}} L_{\Sigma'} [ x ]$. Now minimise $\Sigma' 
$to a $\Sigma > \zeta$ with $L_{\zeta} [ x ] \prec_{\Sigma_{2}} L_{\Sigma} [ x
]$, thus $( \zeta , \Sigma )$ is as required.

Conversely: if it is the case that $\pe ( x ) \da$ the conclusion is trivial,
so suppose otherwise and that $( \zeta , \Sigma )$ is some
$x$-$\Sigma_{2}$-extendible pair satisfying the right hand side. By
$\Sigma_{2}$-extendibility, $\Lambda ( e,x, \Sigma )$ is also $0$. By the
$\liminf$ rule the snapshot of $\pe ( x )$ - which we can envisage running
inside $L_{\Sigma} [ x ]$ - at time $\zeta$ is $\Sigma^{L_{\zeta} [ x ]}_{2} (
x )$. Again by $\Sigma_{2}$-extendibility, it is the same at time $\Sigma$.
Notice that any cell of the tape, $C_{i}$ say, that changes its value even
once in the interval $( \zeta , \Sigma )$, will, by $\Sigma_{2}$-reflection,
do so unboundedly in both $\zeta$ and $\Sigma$. Consequently we have final
looping behaviour in the interval $[ \zeta , \Sigma ]$. Hence we have our
criterion for `$\ensuremath{\operatorname{eJ}}$-convergence'.
{\hspace*{\fill}}{\qed}

\begin{lemma}
  Suppose we have a 2-nesting $\zeta_{0} <  \zeta_{1} < \Sigma_{1} < 
  \Sigma_{0}$. Suppose at time $\zeta_{0}$ \ of the absolute computation of
  $\p{e} ( m )  $either $\p{e} ( m )$ or a subcomputation thereof, is not yet
  convergent and is at level k of its computation tree. Then at time
  $\zeta_{1}$ it is not yet convergent and control is at a level $\geq$ $k+1$.
\end{lemma}

{\pf} Suppose $k=0$. By $\Sigma_{2}$-reflection and the $\liminf$ rule, $\p{e}
( m )  $is still running, and control is still at depth $k$ at $\Sigma_{0}$.
This mean the snapshots at $\zeta_{0}$ and $\Sigma_{0}$ are identical and thus
$\p{e} ( m )$ has its first loop at $( \zeta_{0} , \Sigma_{0} )$, and the
computation is convergent, and is then effectively over. Suppose for a
contradiction that control is at level 0 also at $\zeta_{1}$ (and again also
at $\Sigma_{1}$). So again $\p{e} ( m )$ has looping snapshots at $( \zeta_{1}
, \Sigma_{1} )$. However this is a $\Sigma_{1}$-fact about $\p{e} ( m )$ that
$L_{\Sigma_{0}}$ sees: ``{{\em There exists a 2-extendible pair $( \bar{\zeta}
, \bar{\Sigma} \nobracket$\/}}) {{\em with $\p{e} ( m )$ having identical
snaphots at level 0 at\/}} $( \bar{\zeta} , \bar{\Sigma} )$.'' \ But then
there is such a pair $\bar{\zeta} < \bar{\Sigma} < \Sigma_{0}$ and {\p{e}}$( m
)$'s computation is again convergent at $\bar{\Sigma}$ contrary to assumption.

The argument for $k \geq 1$ is very similar: if $\liminf_{\alpha \rightarrow
\zeta_{0}} \Lambda ( e,m, \alpha ) = \Lambda ( e,m, \zeta_{0} ) =k$, then
$\liminf_{\alpha \rightarrow \Sigma_{0}} \Lambda \left( \p{e} ( m ) , \alpha
\right) =k$ also. Again, if it entered the interval $( \zeta_{1} , \Sigma_{1}
)$ at this same level $k$ it would loop there, and by the same reflection
argument applied repeatedly would do so not just once but unboundedly below
$\zeta_{0}$ at the same level $k$. But after each successful loop at level
$k$, control passes up to level $k-1$. However then $\liminf_{\alpha
\rightarrow \zeta_{0}} \Lambda ( e,m, \alpha ) =k-1$. Contradiction!
{\hspace*{\fill}}{\qed}

\begin{lemma}
  \label{Boundedness}{{\em (Boundedness Lemma for computations recursive in
  $\tmop{eJ}$)\/}} Let $\beta_{0}$ be the least infinitely nested ordinal in
  some ill-founded model $M$ with $\ensuremath{\operatorname{WFP}} ( M )
  =L_{\beta_{0}}$. Let $\alpha_{0}$ be least with $L_{\alpha_{0}}
  \prec_{\Sigma_{1}} L_{\beta_{0}}$. Then any computation
  $P^{\ensuremath{\operatorname{eJ}}}_{e} ( m )$ which is not convergent by
  time $\alpha_{0}$, is divergent.
\end{lemma}

{\nod}{\pf} Let $\zeta_{0} < \cdots  <  \zeta_{n} <  \cdots   \, \beta_{0} \, 
\cdots   \subset  s_{n} \subset \cdots \subset  s_{0}$ witness the infinite
nesting at $\beta_{0}$ in $M$. By the definition of $\alpha_{0}$ no $\p{e} ( m
)$ is convergent at a time $\alpha \in [ \alpha_{0} , \beta_{0} )$ as this
would be a $\Sigma_{1}$-fact true in $L_{\beta_{0}}$; but then by
$\Sigma_{1}$-reflection, it is true in $L_{\alpha_{0}}$. But if $\p{e} ( m )$
is not divergent before $\beta_{0}$, it will be by $\beta_{0}$: the previous
lemma shows that $\Lambda ( e,m, \zeta_{n} ) < \Lambda ( e,m, \zeta_{n+1} )$ \
holds in $M$. \ But these level facts are absolute to $V$, as they are
grounded just on the part of the absolute computation tree being built in
$L_{\beta_{0}}$ as time goes towards $\beta_{0}$ (and are not dependent on
oracle information from $\ensuremath{\operatorname{eJ}}^{M}$ which perforce
will differ from the true $\ensuremath{\operatorname{eJ}}$); so $\p{e} ( m
)$'s computation tree will have an illfounded branch at time $\beta_{0}$.
{\hspace*{\fill}}{\qed}

The above then shows that the initial segment $L_{\alpha_{0}}$ of the
$L$-hierarchy contains all the information concerning looping or convergence
of computations of the form $P^{\ensuremath{\operatorname{eJ}}}_{e} ( m )$. A
computation may then continue through the wellfounded part of the computation
tree for the times $\beta < \beta_{0}$ but if so, it will be divergent.
Relativisations to real inputs $\vec{x}$ are then straightforward by defining
$\beta_{0} ( \vec{x} )$ as the least such that there is an infinite nesting
based at that ordinal in the $L_{} [ \vec{x} ]$ hierarchy {\etc}

\begin{lemma}
  Let $x \sset \omega$. Then $T_{\Sigma ( x )}^{2} ( x ) \dfs
  \Sigma_{2}$-$\ensuremath{\operatorname{Th}} ( L_{\Sigma ( x )} [ x ] )$ is
  $\ensuremath{\operatorname{eJ}}$-recursive in $x$.
\end{lemma}

{\nod}Proof: There is an index $e$ so that running $P_{e} ( x )$ asks in turn
if ?$n \in T_{\Sigma ( x )}^{2} ( x ) ?$ for each $n$, and will receive a
$0/1$ answer from the oracle $\ensuremath{\operatorname{eJ}}$. Consequently
$P_{e}$ may compute this theory on its output tape, and then halt.
{\hspace*{\fill}}{\qed}

{\nod}Remark: $T_{\Sigma ( x )}^{2} ( x ) \equiv_{1 } \tilde{x}$ (by Fact 2
(iii) above).

\begin{lemma}
  \label{model1}Let $x \sset \omega$. Then a code for $L_{\Sigma ( x )} [ x ]$
  is $\ensuremath{\operatorname{eJ}}$-recursive in $x$.
\end{lemma}

{\nod}Proof: There is a standard ittm program that on input $\tilde{x}$ will
halt after writing as output a code for $L_{\Sigma ( x )} [ x ]$. Thus, by the
last remark and lemma, a code for $L_{\Sigma ( x )} [ x ]$ is also
$\ensuremath{\operatorname{eJ}}$-recursive in $x$.{\hspace*{\fill}}{\qed}

Further:

{\bu} \ (i) For any $e,x$, the first repeating snapshot $s ( e,x )$ of $P_{e}
( x )$ is $\ensuremath{\operatorname{eJ}}$-computable in $x$, as is a code for
$L_{\rho_{0}} [ x ]$, $L_{\rho_{1}} [ x ]$ and $L_{\rho_{1}^{+}} [ x ]$ where
$\rho_{0} , \rho_{1}$ are the ordinal stages of appearance of the first
repeating snapshot $s ( e,x )$, and $\rho_{1}^{+}$ is the least $\bar{\rho} >
\rho_{1}$ which is a limit of $s ( e,x )$-admissibles.

{\bu} We may thus have subroutines that ask for, and compute such objects
during the computation of some $\p{f} ( y )$ say. Since satisfaction is also
ittm-computable, we may query simply whether ?$L_{\rho_{1}} [ x ] \models
\sigma$? and receive an answer.

One may show:

\begin{theorem}
  \label{equiv}Any two of the functionals $E 
  ,\ensuremath{\operatorname{eJ}}$, and $\ensuremath{\operatorname{iJ}}$ are
  mutually ittm-recursive in each other.
\end{theorem}

{\pf}This uses, in the direction to obtain $\ensuremath{\operatorname{iJ}}$ or
$\ensuremath{\operatorname{eJ}}$ recursive in $E$, an appropriate version of
the Normal Form Theorem from {\cite{W09}}. {\hspace*{\fill}}{\qed}

We collect together some of the above Facts and results, in order to
abbreviate our descriptions of algorithms This will help to have a library of
basic algorithms which we shall simply quote as being `recursive in
$\ensuremath{\operatorname{eJ}}$' without further justification.

\begin{definition}
  (Basic Computations-$\ensuremath{\operatorname{BC}}$) (i) Any standard
  ittm-computation $P_{e} ( n,x )$ is Basic.
  
  (ii) If a code for an $\alpha  $ordinal is given, then the computations that
  compute: a) \ for any $x$ (a code for) $L_{\alpha} [ x ]$ b) the
  satisfaction relation for $L_{\alpha} [ x ]$ is Basic (in (the code for)
  $\alpha$); (and shows those objects are
  $\ensuremath{\operatorname{eJ}}$-recursive, if $\alpha$ is).
  
  The following are all $\ensuremath{\operatorname{eJ}}$-recursive, and Basic:
  
  (iii) The function $x \rightarrowtail \tilde{x}  $;
  
  (iv) The function that computes $x \rightarrowtail \Sigma ( x )$, the larger
  of the next extendible pair in $x;$
  
  (v) \ The function that computes $x \rightarrowtail \Sigma ( x )^{+}$;
  
  (vi) Any others that we may need to add.
\end{definition}

Stronger ordinals than simply $\Sigma ( x )^{+}$ can be
$\ensuremath{\operatorname{eJ}}$-recursive:

\begin{lemma}
  \label{4.30} There is a recursive sequence of indices $\pa{e_{i} | 0 \leq i<
  \omega }$ so that for any $\alpha < \omega_{1}$ with a code $x \in \cant$,
  $P^{\ensuremath{\operatorname{eJ}}}_{e_{i}} ( x )$ computes a code for the
  next $i$-extendible $\zx{i} > \alpha$. 
\end{lemma}

{\pf} For $i=0$ this has been done using Basic Computations. Suppose $e_{i}$
has been defined, and we describe the programme
{\p{e\ensuremath{_{\textrm{i+1}}}}}. Assume without loss of generality that
$\alpha =0$, $x=\ensuremath{\operatorname{const}}_{0}$. Then $\p{e_{i}} ( 0 )$
computes a code for the least $i$-extendible, $\zeta_{0} :=${\zx{i}} say. By a
basic computation let a slice of the scratch tape $R$ be designated to hold
$T_{\zeta_{0}}^{2}$; $R:=T^{2}_{\zeta_{0}}$. A code for $\zeta_{0}$ is
recursive in $T^{2}_{\zeta_{0}}$. Now compute $\p{e_{i}} ( R )$. This yields
the next $i$-extendible $\zeta_{1} = ^{i} \zeta_{1}$. Now, using Basic
Computations, write successively to $R$ the theories $T^{2}_{\zeta_{0}}
,T_{\zeta_{0} +1}^{2}   , \ldots  ,T^{2}_{\zeta_{0} + \beta}    , \ldots$ for
$\beta < \zeta_{1}$. We note that at limit stages $\lambda \leq \zeta_{1}$,
$R$ will contain ``liminf'' theories $\hat{T}_{\lambda}
=\ensuremath{\operatorname{Liminf}}_{\alpha \rightarrow \lambda}
T^{2}_{\alpha}$ (by the usual automatic ittm liminf process) but that
$T^{2}_{\lambda}$ is uniformly r.e. in \ $\hat{T}_{\lambda}$. (For the latter
see Fact 2. It is easy to argue that $\hat{T}_{\lambda} \supseteq
T^{2}_{\lambda}$, and that if $\ensuremath{\operatorname{supS}}^{1}_{\lambda}
= \lambda$ then we have equality, it is the bounded case of $S^{1}_{\lambda}$
in $\lambda$ that requires argument. The point of this exercise of writing
theories to $R$ is to ensure continuability of the computation, and that we do
not start to loop too early. The `writing out' of all levels of the theories
to $R$, is a precautionary step: in general we do not have $\hat{T}_{\zx{i+1}}
= \liminf_{^{i} \zeta \rightarrow \zx{i+1}  } \hat{T}_{^{i} \zeta}$.) And
again a code for $\lambda$ is then recursive in $T^{2}_{\lambda}$.

Set $R :=$ $\hat{T}_{\zeta_{1}}$; by the comments just made
$T^{2}_{\zeta_{1}}$ is r.e. in $R$ and $R \in L_{\zeta_{1} +1}$ (this is why
we are writing out these theories, to ensure that we loop at our desired
target); now compute $\p{e_{i}} ( R )$ and repeat this process. \ As there is
no means for the machine to halt, there is a least looping pair $( \zeta ,
\Sigma )$ of ordinals. \ Let $\left( \zx{i+1} \right.$,$\left. \sx{i+1}
\right)$ be the least $i+1$-extendible pair. \ We claim that this is the pair
$( \zeta , \Sigma )$. Suppose $\zeta < \zx{i+1}$. By the repetition of the
contents of $R$ in the loop points, we have $\hat{T}_{\zeta} =
\hat{T}_{\Sigma}$, in the above algorithm, hence $T^{2}_{\zeta}
=T^{2}_{\Sigma}$, \ and thus $L_{\zeta} \prec_{\Sigma_{2}} L_{\Sigma}$. But
then $\zeta$ is an extendible limit of $i$-extendibles, as $\zeta$ is a limit
point of this looping process. This contradicts the minimality of $\zx{i+1}$.
Hence $\zeta$ equals the latter, and $\Sigma = \sx{i+1}$ follows.

Hence we may compute $\hat{T}_{\zx{i+1}}$, \ $\zx{i+1}$ by means of an
eventually stabilizing looping programme. We let
{\p{e\ensuremath{_{\textrm{i+1}}}}} be the programme just described followed
by the basic comp. that finds a code for $\zx{i+1}$ by a method uniformly r.e.
in $\hat{T}_{\zx{i+1}}$.

Finally note that the continuing description of the programme $\p{e_{i+2}}$
from {\p{e\ensuremath{_{\textrm{i+1}}}}} merely repeats the above but altering
only a few indices. We may thus determine a recursive function $i \mapsto
e_{i+1}$. {\hspace*{\fill}}{\qed}

Entirely similar is:

\begin{lemma}
  \label{nextext}There is a (Turing) recursive sequence of indices $\pa{e'_{i}
  \mid i< \omega}$ so that $\p{e'_{i}}( x )$ writes a code for $^{i}
  \Sigma^{} ( x )$, the least $\Sigma_{2}$-extension of $L_{ ^{i} \zeta} [ x
  ]_{}$. 
\end{lemma}

\section{The determinacy results}

We shall assume a certain amount of familiarity of working with ittm's and
shortcuts amounting to certain subroutines, so as not to overload the reader
with details.

\begin{theorem}
  \label{overall}For any $\Sigma^{0}_{3}$ game $G ( A;T ) ,$(with $T$ say
  recursive) if player {\pone} has a winning strategy, then there is such a
  strategy recursive in $\ensuremath{\operatorname{eJ}}$; if player {\ptwo}
  has a winning strategy, then there is such a strategy either recursive in
  $\ensuremath{\operatorname{eJ}}$, or else definable over $L_{\beta_{0}}$.
\end{theorem}

{\nod}{\textbf{Proof of \ref{overall}}}

Idea: We suppose $A= \bigcup_{n} B_{n}$ with each $B_{n} \in \Pi^{0}_{2}$,
with an initial game tree $T$. \ For expository purposes we shall assume that
$T= ^{< \omega} \omega$ - relativisations will be straightforward. We shall
provide an outline of a procedure which is recursive in
$\ensuremath{\operatorname{eJ}}$ and which will either provide a strategy for
{\pone} in $G ( A;T )$ (if such exists) or else will diverge in the attempt to
find a strategy for {\ptwo}. We wish to apply the main Lemma 3 of
{\cite{We11}} for the successive $B_{n}$. The control of the procedure will be
at different {{\em Levels\/}} of the initial finite path tree of the
computation. At Level 0 will be the main process, but also the procedure for
finding witnesses and strategies involved in the arguments for the Main Lemma
applied with $B=B_{0}$. We first search for a level in the $L$-hierarchy whose
code is $\ensuremath{\operatorname{eJ}}$-recursive~and for which we can define
a non-losing subtree $T' \sset T$, for which all $p \in T'$ have witnesses
$\hat{T}_{p}$ to $p$'s goodness in the sense of (i) and (ii) above. In fact we
shall search for pairs of levels in the $L$-hierarchy, in the sequel, between
which we have absoluteness of our non-losing subtrees. After having found
such, this data will be encoded as a real (these routine details, the reader
will be pleased to learn, we omit) and a subroutine call made to a process at
the {{\em lower\/}} {{\em Level 1\/}} which will attempt to find the right
witnesses {\etc} to apply the Lemma for $B=B_{1}$. We now search for a further
level of the $L$-hierarchy which again has the right witnesses to goodness to
all the possible relevant subtrees associated with positions $p_{2}$ of length
$2$. \ As we search for such an $L_{\alpha}$, we may find that some of our
original witnesses to goodness at Level 0 no longer work in our new
$L_{\alpha}$, or even more simply that our $T'$ from Level 0 now has nodes $p$
which have become winning for {\pone} in this $L_{\alpha}$. We accordingly
keep testing the data handed down to see if any of it has become `faulty' in
this respect. If so, then we throw away everything we have done at Level 1,
but pass control back up to Level 0 together with the ordinal height of the
current $L_{\alpha}$ we reached. We then go back to searching for an
$L_{\alpha'}$ which is `good' in all of these previous respects at Level 0 for
a new $T'$, which we then shall pass down to Level 1 for another attempt.

Eventually we shall reach a stage where we have a sufficiently large model
where all the data and our witnessing subtrees work at both levels 0 and 1.
Accordingly again all {{\em this\/}} data is passed down to the subroutine at
{{\em Level 2\/}} for assessing potential subtrees for application in the
Lemma to be applied for $B=B_{2}$. Proceeding in this fashion, testing as we
go the validity of our data trees {{\em en route\/}} and passing back up to
the Level of the tree that has failed if so, we find we work at increasing
depth - that is at lower Levels $n$ with increasing $n$. \ If {\ptwo} has a
winning strategy then there will be an infinite path descending through all
the Levels and hence the computation will diverge. One point will be to remark
that if $I$ has a winning strategy then this process will discover it: this
requires us checking that we don't come up against a `wall' in the ordinals
$\alpha$ so that we cannot find a code for an ordering of a longer order type
- because our computation has stabilized, or in other words is in a loop, and
we are stuck below the length of that loop.

Hence if there is no such wall, and $G ( A;T )$ has a winning strategy for
{\ptwo} only definable over $L_{\beta_{0}}$. then we can theoretically keep
computing ordinals up to $\beta_{0}$.

Our task now is to achieve a balance between{\color{blue} } giving enough of
these details that the reader is convinced, and without causing the eye to
glaze over with overwhelming (and unnecessary) {{\em minutiae\/}}.

{{\em In general:\/}} given a tree $S$ in a model $M$, used in a game $G (
\bar{A} ,S ) $, and without a strategy for player \pone in $M$, then we
shall denote the subtree of non-losing positions for {\ptwo} in $M$ by
$S'^{M}$ (or just $S'$). For $R \in \power \left( \nat \right)$, $\tau^{+} ( R
)$ will denote the sup of the first $\omega$ many $R$-admissibles beyond
$\tau$. By $\Sigma^{k} ( R )$ we shall mean, where $\zeta^{k} ( R )$ is the
least $k$-extendible in the $L_{\alpha} [ R ]$ hierarchy, that $\Sigma^{k} ( R
)$ is the least ordinal with $L_{\zeta^{k} ( R )} [ R ] \prec_{\Sigma_{2}}
L_{\Sigma^{k} ( R )} [ R ]$. If $k=1$ we drop it and write simply $\zeta ( R
)$ {\etc}We note that if {\color{blue} {\color{blue} {\color{blue}
{\color{blue} {\color{blue} {\color{blue} }}}}}}$L_{\Sigma^{+} ( R )} [ R ]$
has no proper $\Sigma_{1}$-substructures, then $T_{\Sigma^{+} ( R )}^{1} [ R
]_{}$ {\dfs} $\Sigma_{1}$-Th($L_{\Sigma^{+} ( R )} [ R ]$) - in the language
of set theory with a predicate symbol for $R$ - is not in $L_{\Sigma^{+} ( R
)} [ R ]$; moreover (ordinarily) recursive in $T_{\Sigma^{+} ( R )}^{1} [ R ]$
is a wellorder of type $\Sigma^{+} ( R )$. We shall let the notation
$^{\alpha} M$ vary over structures of the form $L_{^{\alpha} \Sigma^{+}} [ T
]$.

[Commentary is provided in square brackets following a \% sign.]

As a warm-up we prove the following lemma using just Basic Computations.

\begin{lemma}
  \label{4.29}There is a computation that on input codes for
  $T,${\bfseries{{\pa{B\tmrsub{n}}}}} will halt either with a winning strategy
  for \pone, or else with an encoded $T'$ - the set of non-losing positions
  for {\ptwo} in $G ( A;T )$- membership of which is absolute between some
  $L_{\zeta} [ T ]$ and $L_{\Sigma^{+}} [ T ]$.
\end{lemma}

(0): {{\em  We commence with cutting up recursive infinite disjoint slices of
the scratch tape to be reserved as `registers' for the reals coding
{\pa{B\tmrsub{n}}}, $T,T'$,$\Sigma^{+}$,$\ldots$ \/}}{\color{red} }, (and more
such will be needed at lower Levels, as data is passed down in the argument
that follows, but we shall not mention these, rather leave it to the reader to
do the preparatory mental scissor work).\\

{\bu} {{\em Set:\/}} \ $T' :=T$.

\nod (1)  {\bu} \ {{\em Compute: \/}}$M:=L_{\Sigma^{+} ( T' )} [ T' ]$,{{\em  \/}}
and set \ $\Sigma^{+} := \Sigma^{+} ( T' )$.

{\bu} ?$T'^{M} = \emp ?${{\em  If $T'^{M} = \emp$ then 
\pone \, has a winning strategy in $G ( A;T ) ^{M}$, and
this may be found in $M$ and printed out on the output tape; then \/}}STOP.
Otherwise CONTINUE.\\

 [\% As $M$ is a model of $\ensuremath{\operatorname{KPI}}$ such a
winning strategy is winning in $V$.] \\

{\bu} {{\em \/}}? {{\em Is\/}} $T'^{M} = T'  ?$ \ \

\nod (2) {\bu} {{\em If\/}} NO {{\em then $T' \supset T'^{M}$\/}} {{\em and then
some winning strategies are newly available to \pone in $M$ that are for
some $p \in T' \backslash T'^{M}$. \/}} {{\em Set\/}} $T' := T'^{M}  $; \ \
GOTO (1).\\

[\% Note that the new $T'$ is a proper
subtree of the old.]\\

{\bu} {{\em If\/}} YES, {{\em then we may\/}} STOP {{\em with a
suitable $T'$ encoded in its register.\/}}\\

[\% \ Of course in order to obtain $M, \Sigma^{+}  ,${\etc} this officially
requires a call to a subcomputation at the next level down, but the above is
just a schematic description of the process, and so we suppress that level of
detail. The point is that the $T'$ are a decreasing sequence of sets. Hence
keeping track of these $T'$ at the top level suffices for the procedure to
continue: we don't need to keep track of, {\eg}, the ordinals heights of the
structures $M$, and the concomitant worries about the liminf action at limit
stages. Thus the above can be all effected using Basic Computations (and
variants thereon).]\\

{{\em Claim 1\/}} {{\em Either the program halts with a winning strategy for
{\pone} in $G ( A;T )$ or, at some point strictly before the next
$2$-extendible above $\Sigma^{+} ( T )$ in the cycle, the answer to the query
\ {\tmem{? Is}} $T'^{M} =T'  ?$ is affirmative.\/}}

{\pf} Note first that the computation uses only BC's and each of these only
require a computation of length the next extendible pair at most. Suppose $(
\zeta_{0} , \Sigma_{0} )$ is any extendible pair that is a limit of such,
above $\Sigma^{+} ( T' )$. \ We imagine the computation as being performed as
a $\Sigma_{2}$-recursion in $T$ in $L_{\Sigma_{0}} .$ Then suppose, for a
contradiction, that by the $\zeta_{0}$'th turn through the cycle, we have not
had an affirmative answer. In the $\nu$'th turn through the cycle (for $\nu <
\zeta_{0}$) let $T'$ be denoted by $T'_{\nu}$. Then the $T'_{\nu}$, as
remarked, are strictly decreasing. \ Now by an easy reflection argument, one
sees that on a tail of $\nu < \zeta_{0}$, the $T'_{\nu}$ must be the same. [If
``$\all \nu \ex \nu' > \nu \ex p \left( p \in T'_{\nu'} \back T'_{\nu' +1}
\right)$'' holds in $L_{\zeta_{0}}$ it will also hold in $L_{\Sigma_{0}}$. But
if $p_{0} \in T'_{\nu'} \back T'_{\nu' +1}$ the $\nu'$ for which that happens
is $\Sigma_{2}$-definable in $L_{\Sigma_{0}}$ from $p_{0}$; but that implies
$\nu' < \zeta_{0}$. This contradicts the quoted formula.] So an affirmative
answer must have occurred. \hfill{\qed  \mbox{ }Claim 1 and Lemma.}\\

We now assemble these building blocks to form a programme based on the
argument of the proof of Theorem \ref{gamma=beta} surveyed above.\\

{\nod}{\textbf{Proof of Theorem \ref{overall}:}}

We outline the argument at the various levels of computation in the oracle
calls of a master computation at level $\Lambda =0$. We proceed by describing
the actions of the programmes being called, which the reader may reformulate
as official queries to the $\ensuremath{\operatorname{eJ}}$-functional as
oracle. At the end of the description we justify the claim that this is a bona
fide $\ensuremath{\operatorname{eJ}}$-recursion.

(1) $\Lambda =0$.

{\bu} The master or control programme computes successively lengthening
structures \ $^{1} M=L_{\Sigma^{+}} [ T' ]$ until $T'$ is seen to stabilize
between one such structure $^{1} M$ and the next, $^{1} M'$.

[\% This we saw done effectively by a machine in the proof of Lemma
\ref{4.29}, with $T'$ so stabilizing before the next $2$-extendible. This
process involved oracle queries to Level $\Lambda =1$, but again we suppress
these details.]

{\bu} With $T'$ stabilized, the programme asks the following - when suitably
formulated - oracle query of $\ensuremath{\operatorname{eJ}}$. The query
sub-computation we view as enacted at $\Lambda =1$. We suppose that it is the
computation {\p{e\ensuremath{_{\textrm{0}}}}}$( x )$ where $x= \pa{1,
\pa{B_{n}}_{n} , \tmt{0} , ^{1} M}$ (suitably coded), whose action is
described below starting at (2).\\

$Q^{1} :$ ? {{\em Defining $\tmt{1}$ from the current $T'$ in $\tmt{0}$, do
all the trees {\nod}in $\mathbbm{T}^{1}$ become eventually settled \/}}?\\

\nod [\% Recall that the trees of {\tmt{\textsuperscript{1}}} are of the form:

a) \ \ $\hat{T}_{p}$ \ \ \ \ ({\dfs} the current $^{1} M$-least witness to the
goodness of $p \in T'$) and

b) $( \hat{T}_{p} )'  $ ({\dfs} its tree of non-losing positions for {\ptwo});
as well as (where $T^{\ast} \left( \emp \right)$ is set to $\hat{T}_{\emp}$)

c) \ $T^{\ast} \left( \emp \right)_{p_{2}}$ and $\left( \left( T^{\ast} \left(
\emp \right) \right)_{p_{2}} \right)'$ for relevant $p_{2}$.

We adopt the convention, that `` $\mathbbm{T}^{l}$ becomes eventually settled
'' or `` $\mathbbm{T}^{l}$ is stable up to ordinal $\tau$'' to be a shorthand
affirming that all the constituent trees of the family$ \mathbbm{T}^{l}$ are
stable {{\em per\/}} their definitions up to $\tau$.

Note also: that since $T'  $ has survived intact from one $^{1} M$ structure
to the next $^{1} M'$ say, we can deploy the `survival argument' of Lemma
\ref{2.12}; this means that both structures see that all $p \in T'$ are good,
and this is a sufficient criterion for the definition of $\tmt{1}$ over $^{1}
M$ to instantiate all the needed trees, which then exist in $^{1} M$ (indeed
$( \mathbbm{T}^{1} )^{^{1} M} \sset ( L_{^{1} \zeta} )^{^{1} M}$). Hence the
query is therefore immediately meaningful. ]

(2) {\p{e\ensuremath{_{\textrm{0}}}}}$( x )$ answers the query by first taking
from$ x$ the current data, and on seeing the initial flag $1$, computes
successive models $^{1} M$, and keeps a register of the successive theories
$T^{2}_{\alpha}$, of increasing ordinal height in the manner of the proof of
Lemma \ref{4.30}. These operations are using our BC's.\\

If (Case 0): An $^{1} M$ is reached that contains a winning strategy $\sigma$
for {\pone} in $G ( A;T )$ then the programme HALTS and passes $x' =
\pa{\sigma}$ back up to the master programme at $\Lambda =0$;

If (Case 1): $T'$ changes from one structure $^{1} M$ to the next (``$T'$
becomes unstable'') then the programme HALTS and with the current $\tmt{0} =
\pa{T,T'^{ ^{1} \!\! M}}$, passes the current $x' = \pa{1, \pa{B_{n}}_{n} ,
\tmt{0} , ^{1} M}$ back up to the master programme at $\Lambda =0$; and
RETURNS TO (1);

If (Case 2): $T'$ remains stable but some $S$ $\in \mathbbm{T}^{1}  $does not
by the end of the eventual loop in {\p{e\ensuremath{_{\textrm{0}}}}}$( x )$,
then the answer to $Q^{1}$ is ``No'' (or ``$0$'') and $x' = \pa{0}$ and
control are passed back up to the master programme at $\Lambda =0$.\\

In Case 0, the Master programme halts with this $\sigma$ as output.

\ \ \ \ \ \ [\% note that as $M$ is closed under admissibles, $\sigma$ is a
w.s. for {\pone} in $V$.]

In Case 1, the Master programme continues to calculate successive models,
re-starting from the $M$ passed up in $x'$.

In Case 2, the Master programme, on receiving ``No'', and using BC's, computes
the length of the loop just passed, call it $\Sigma$, and then continues
calculating successive models, with the first such in this series containing
the ordinal $\Sigma$.\\

[\% \ Note that: (A) $T'$ must become eventually settled under the repeated
calculation of longer $M$'s by the time of the next (or indeed any) larger
element $\zx{2} \in E^{2}$, or $\zx{\alpha}$ ($ \alpha \geq 2 )$ for that
matter. Hence the loop (1) $\imp$(2) (Case 1) $\imp$(1) will be broken out of
by the time the length of the models $M$ approaches the next $\zx{2}$.

(B) For Case 2: we cannot immediately deploy a shrinking argument on the trees
to conclude that we have stability of all trees in $\mathbbm{T}^{1}$ by the
next extendible, since the actual underlying trees \ $\hat{T}_{p}$, $T^{\ast}
\left( \emp \right)_{p_{2}}$ may be changing. However the eventual loop
$\Sigma$ whose length the Master programme computes, is that of a 2-extendible
in $E^{2}$; this is ensured by the writing out of the theories
$T^{2}_{\alpha}$ in the manner of the argument of the proof of Lemma
\ref{4.30}. If the loop (1) $\imp$(2) (Case 2) $\imp$(1) repeatedly occurs
from some point on, then for all sufficiently large $^{2} \Sigma$ below the
next$ ^{3} \zeta$ (and so also by $\Sigma_{2}$-reflection, below the $^{3}
\Sigma $ corresponding to $\zx{3}$).There is $\mathbbm{T}^{1}$ so that (for
all sufficiently large $^{1} \Sigma < ^{2} \Sigma )  $( $\mathbbm{T}^{1} = (
\mathbbm{T}^{1} )^{L_{^{1} \Sigma}}$) and so we shall end up in Case 3 below.
] \\

The last possibility is:\\

If (Case 3): All \ $S \in \mathbbm{T}^{1}  $ become stable between two successive
$^{1} M$-structures, $M_{1} ,M_{2}$.

{\nod}The sub-computation now makes in turn a further query sub-computation
which in turn we view as enacted at $\Lambda =2$. We suppose that it is the
computation {\p{e\ensuremath{_{\textrm{0}}}}}$( x )$ where we collect the
current values
$$\tmt{1} =  \langle \pa{\hat{T}_{p} \mid p \in
T'} , \pa{( \hat{T}_{p} )' \mid p \in T'},\langle T^{\ast} \left( \emp
\right)_{p_{2}} , \left( T^{\ast} \left( \emp \right)_{p_{2}} \right)' \mid 
p_{2}\mbox{ relevant }\rangle \rangle$$
{\nod}and set:
$$x=  \pa{2, \pa{B_{n}}_{n} , \tmt{0} ,
\tmt{1} ,M_{1}};$$
{\nod}and where the query is:\\

$Q^{2} :$ ? {{\em Defining {\tmt{2}} from the current $\tmt{0} , \tmt{1}$ of
$x$, \ do all the trees $S$ in $\mathbbm{T}^{2}${\nod} become eventually
settled\/}} \ ?\\

[\% \ Just as following $Q^{1}$, the stability of all the trees $( \hat{T}_{p}
)'$ and $\left( T^{\ast} \left( \emp \right)_{p_{2}} \right)'$ from one model
to the next guarantees the existence of all the trees of $\tmt{2}$ by the
survival argument. ]\\

(3) {\p{e\ensuremath{_{\textrm{0}}}}}$( x )$ is programmed so that when it
takes from$ x$ the current data, and sees the initial flag $2$, it will
continue to compute successive models $^{2} M$, (which it can by Lemma
\ref{nextext}) and write out theories as before, using Basic Comps, but now
act as follows.\\

If (Case 0): \ A $^{2} M$ contains a winning strategy $\sigma$ for {\pone} in
$G ( A;T )$ then this sub-computation HALTS and passes $x' = \pa{\sigma}$ back
up to the programme at $\Lambda =1$;

If (Case 1): $T'$ becomes unstable, then the subcomputation HALTS and passes
the current $x' = \pa{0,T,T' , ^{2} M}$ back up to the programme at $\Lambda
=1$;

If (Case 2): $T'$ remains stable but some $S \in \mathbbm{T}^{1}$ does not at
some stage, between two successive models $^{2} M_{1} , ^{2} M_{2}$, then the
subcomputation HALTS and the current $x' = \pa{2, \pa{B_{n}}_{n} , \left(
\tmt{0} , \tmt{1} \right)^{^{2} M_{2}} , ^{2} M_{2}}$ with the current values
of the data, and control, are passed back up to $\Lambda =1$;

If (Case 3): $T'$ and all $S \in \mathbbm{T}^{1}$ remain stable but some $S
\in \mathbbm{T}^{2}$ do not, then the answer to $Q^{2}$ is ``No''.\\

In Cases 0,1 the relevant information will be passed up in turn to the master
computation at $\Lambda =0$ and will be acted on appropriately.

In Case 2, the sub-computation at $\Lambda =1$, restarts using BC's, and
computes structures $^{1} M$ as at (2).

In Case 3, the sub-computation at $\Lambda =1$, is programmed to use BC's, to
compute the length of the loop just passed, say to $\Sigma$, and then
continues calculating successive models in the usual manner as at (2), with
the first such in this series containing the ordinal $\Sigma$.

[\% \ Note that: the comments on the loops at (A), (B) will hold here.
Additionally:

(C) If the loop (2) $\imp$(3) (Case 3) $\imp$(2) occurs from some point on,
then for sufficiently large $^{3} \Sigma$ below the next $ ^{4} \zeta$, (and
so also by $\Sigma_{2}$-reflection, below the corresponding $^{4} \Sigma  )$
there is $\mathbbm{T}^{2}$ so that (for sufficiently large $^{2} \Sigma < ^{3}
\Sigma   ( \mathbbm{T}^{2} =\mathbbm{T}^{2} )^{L_{^{2} \Sigma}}$) and so we
shall end up in Case 4 below. ]\\

The last possibility is:\\

\pagebreak
If (Case 4): {\tmt{2}} becomes stable between two successive $^{2} M$-structures,
$M_{1} ,M_{2}$.\\

Again, {\nod}the current sub-computation makes a query sub-computation which
in turn we view as enacted at $\Lambda =3$. We suppose that it is the
computation {\p{e\ensuremath{_{\textrm{0}}}}}$( x )$ where we set
$$\tmt{2} = \pa{\hat{T} ( p_{2} )_{p} , \hat{T} ( p_{2} )_{p}' \mid p \in
\left( \tptwo \right)'} ,  \pa{T^{\ast} ( p_{2} )_{p_{4}} , ( T^{\ast} ( p_{2}
)_{p_{4}} )' \mid  p_{4}  \ensuremath{\operatorname{relevant}}}$$ and
$$x= \pa{3,
\pa{B_{n}}_{n} , \tmt{0} , \tmt{1} , \tmt{2} ,   M_{2}}$$

\nod and where the query is:\\

$Q^{3} :$ ? {{\em Defining {\tmt{3}} from the current $\tmt{0} , \tmt{1} ,
\tmt{2}$ in $x$, does $\mathbbm{T}^{3}${\nod} become eventually settled\/}} \
?\\

We hope the reader will have seen the pattern emerging in this description of
the programme $P_{e_{0}}$. However the reader is entitled to ask: have we
described a genuine programme for such oracle machines? And secondly, what is
the outcome?\\

\begin{figure}[htbp]
\resizebox{12.2cm}{!}{\input{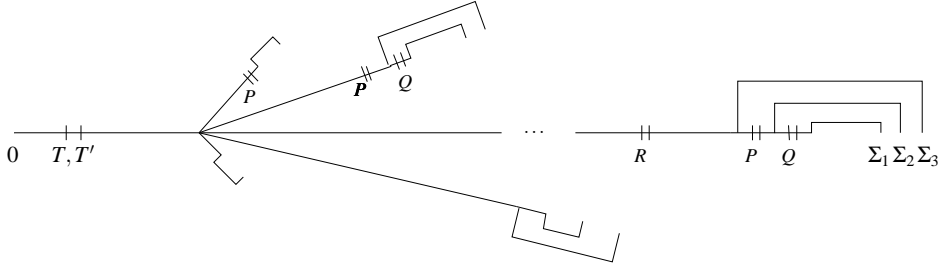}}
 \begin{picture}(0,0)%
\end{picture}%
\setlength{\unitlength}{3947sp}%

\caption{In this diagram $T'$ is stable up to (but not beyond) $\Sigma_{3}$.}
\label{figure:example}
\end{figure}


A typical $3$-nesting diagram is at Figure 1. $T'$ is assumed to be
stable up to $\Sigma_{3}$. Thus beyond the branch given, there are no winning
strategies for {\pone} for any $T'_{p}$ for any $p \in T$ appearing in the
interval beyond the branch point up to $\Sigma_{3}$ (but such may appear in $
L_{\Sigma_{3} +1} )$. Because $T'$ is this long-lived at positions labelled
$P$, we can have all the relevant trees $( \hat{T}_{p} )'$, $T^{\ast} \left(
\emp \right)_{p_{2}}$ and $( ( T^{\ast} \left( \emp \right)
)_{p_{2}} )'$ ({\ie} $\mathbbm{T}^{1}$) \ occurring, and
themselves are stable up to the end of the extendible loop below which they
occur. At the first $2$-nesting illustrated because all the $( (
T^{\ast} \left( \emp ) )_{p_{2}} \right)'$ at $P$ survive to the
end of the outermost nesting, and so {{\em beyond\/}} the top of the inner
nesting, we may conclude that at a position such as $Q$, all the relevant
trees $T^{\ast} ( p_{2} )_{p_{4}} , ( T^{\ast} ( p_{2} )_{p_{4}} )'$ of
$\mathbbm{T}^{2}$ occur below the inner extendible $\zeta$ that starts the
inner nesting loop. The analysis at the $3$-nesting is similar: since $T'$
survives beyond $\Sigma_{2}$, the $\mathbbm{T}^{1}$ trees can be found at
locations $P$; as the $\mathbbm{T}^{1}$ trees, $( \hat{T}_{p} )'$, $T^{\ast}
\left( \emp \right)_{p_{2}}$ and $( ( T^{\ast} ( \emp )
)_{p_{2}} )'$, survive beyond $\Sigma_{1}$, the $\mathbbm{T}^{2}$
trees can be found at locations $Q$. If we had assumed that $T'$ survived
beyond $\Sigma_{1}$ then we could have obtained a shift, with the
$\mathbbm{T}^{1}$ trees obtainable at $R $, the $\mathbbm{T}^{2}$ trees at $P$
and then gone on to find the $\mathbbm{T}^{3}$ trees at $Q.$

We could easily enough have written down $Q^{k+1}$ which, given $\tmt{0} ,
\ldots    , \tmt{k}$ from an $x$ would have formulated definitions for $T^{*}
(p_{2 ( k-1 )} ) \dfs \hat{T} ( p_{2 ( k-1 )} )_{\emp}$, relevant $p_{2k} ,$
and then asked if trees $\hat{T} ( p_{2k} )_{p}$ (being the current $^{k+1}
M$-least witness to the goodness of $ p \in T^{*} ( p_{2 ( k-1 )} )_{p_{2k}}
)'$ and $\hat{T} ( p_{2k} )_{p}'$ (the latter's subtree of non-losing
positions for {\ptwo}), that is the trees of $\mathbbm{T}^{k+1}$, \ became
eventually settled. The required definitions and stability questions are then
entirely uniform in $k$. Hence the instructions for the programme $P_{e_{0}}$
on input an $x$ coding some $\left\langle k+1, \pa{B_{n}}_{n} , \tmt{i}   ( i
\leq k )  ,^{k} M \right\rangle$ may be effectively written down in terms of
$k$ and the given tuple of data. It is enacted by considering successive
$^{k+1} M$ structures, and by writing down theories $T^{2}_{\alpha}$ as
before. The number of Cases to be considered at query $Q^{k+1}  $ is $k+3$:
Cases (0)-($k$) result in a HALT at that level $\Lambda =k+1$, with an
effectively determined $x'$ to be passed up to the level $\Lambda =k$ above;
whilst Case $k+2$ requires returning to $\Lambda =k$ and computing lengths of
loops {\etc} The final Case $k+3$ is the one of eventual interest and triggers
the query $Q^{k+2}$. Each $Q^{k+1}$ is officially a query of the form
$?\ensuremath{\operatorname{eJ}} ( ( e'_{k+1} ,x ) ) =0/1?$ about how the next
subcomputation loops, and we calculate the relevant $x $ from our data. The
instructions that $e_{k+1}'$ codes include of course those for calculating
$e_{k+2}'$ ready for the next query. However we may argue as below, that these
calculations may be assembled into, or considered as, one whole calculation
embodied in one $\varphi^{\ensuremath{\operatorname{eJ}}}_{e_{0}}$.

In the following we let ``$\all^{\ast} \alpha < \Sigma \varphi ( \alpha )$''
abbreviate ``For all sufficiently large $\alpha < \Sigma   \varphi ( \alpha
)$''. We shall say ``$T'$ is stable up to $^{k} \Sigma$'' to mean
``$T'^{L_{\zeta ( M )}} =T'^{L_{\Sigma^{+} ( M )}}$'' for all sufficiently
large structures $M$ with $\Sigma^{+} ( M ) < ^{k} \Sigma$. This can be
equivalently written as ``$\ex U ( \all^{\ast}  \mbox{}^{k-1} \Sigma < ^{k}
\Sigma ) [ U' = ( T' )^{L_{ ^{k-1} \Sigma}} ]$.''

For $0<l<k$ we shall say ``$\mathbbm{T}^{l}$ is stable up to $^{k} \Sigma$''
to mean $$\ex \mathbbm{T}^{l} \in L_{^{k} \Sigma} \left( \mathbbm{T}^{l}
=\mathbbm{T}^{lL_{\Sigma^{+} ( M )}} \right)$$ for all sufficiently large
structures $M$ with $\Sigma^{+} ( M ) < ^{k} \Sigma$, which, as we have
indicated above, of course is taken, by a convention, to be a shorthand
affirming that all the constituent trees of $\mathbbm{T}^{l}$ are stable {{\em
per\/}} their definitions up to $^{k} \Sigma$.

In the above definition of the algorithm we are employing the following
principle:
\pagebreak

{{\em Suppose $T'  $is stable up to some $^{k} \Sigma$, then  \

for all sufficiently large $ ^{k-1} \Sigma < ^{k} \Sigma \nocomma   (
\mathbbm{T}^{1} \nobracket$ is stable up to $ ^{k-1} \Sigma$ \&

for all sufficiently large $ ^{k-2} \Sigma < ^{k-1} \Sigma \nocomma   (
\mathbbm{T}^{2} \nobracket$ is stable up to $ ^{k-2} \Sigma$ \& $\ldots$

$\vdots$

for all sufficiently large $ ^{2} \Sigma < ^{3} \Sigma \nocomma 
,\mathbbm{T}^{k-2}$ is stable up to $ ^{2} \Sigma$ \&

for all sufficiently large $ ^{1} \Sigma < ^{2} \Sigma \nocomma 
,\mathbbm{T}^{k-1}$ exists $\nobracket \nobracket ) \cdots )$''.\/}}

Less perspicuously but more formally we state this as:

\begin{lemma}
  Suppose $T'  $is stable up to some $^{k} \Sigma$, then  \
  
  $\left( \all^{\ast}  \mbox{}^{k-1} \Sigma < ^{k} \Sigma \right) \left( \ex
  \mathbbm{T}^{1} \right) \left(^{} \all^{\ast}  \mbox{}^{k-2} \Sigma < ^{k-1} \Sigma
  \right)$
  
  $\left[ \mathbbm{T}^{1} = ( \mathbbm{T}^{1} )^{L_{ ^{k-2} \Sigma}} \wedge
  \left( \ex \mathbbm{T}^{2} \right) \left(^{} \all^{\ast} \mbox{} ^{k-3} \Sigma <
  ^{k-2} \Sigma \right) \left[ \mathbbm{T}^{2} = ( \mathbbm{T}^{2} )^{L_{
  ^{k-3} \Sigma}} \wedge \left( \ex \mathbbm{T}^{3} \right)^{} ( \ldots   )  
  \cdots \right. \right.$
  
\noindent  $\ldots \left( \ex \mathbbm{T}^{k-2} \right)\!\!\left(\all^{\ast}\mbox{}^{1}
  \Sigma < ^{2} \Sigma \right)  \hspace{-0.2em}\left( \mathbbm{T}^{k-2}
  = ( \mathbbm{T}^{k-2} )^{L_{ ^{1} \Sigma}} \wedge \left( \ex
  \mathbbm{T}^{k-1} \right) \left( \mathbbm{T}^{k-1} = ( \mathbbm{T}^{k-1}
  )^{L_{ ^{1} \Sigma}} ) ] ] \cdots \right] \right]$.
\end{lemma}

{\pf}Formally by induction on $k$, but the reader may convince themselves of a
representative case, say with $k=3$. \ {\hspace*{\fill}}{\qed}

\begin{note}\em
  The Lemma is really the formal counterpart of the description that precedes
  it. Note that the hypothesis here is fulfilled whenever $^{k} \Sigma$
  approaches some $^{k+1} \Sigma$: for sufficiently large $^{k} \Sigma$ below
  $^{k+1} \Sigma$, $T'$ will be stable even beyond $^{k} \Sigma$.
   Also, by the usual $\Sigma_{2}$-reflection arguments, the above principles
  are equivalent to those obtained by replacing any string ``$< ^{l} \Sigma$''
  by ``$< ^{l} \zeta$''.
  As the program runs there will eventually be subcomputation calls to
  arbitrary levels, as it uses various trees for as long as they survive
  fulfilling their role. But only after $\alpha_{0}$ stages will we be certain
  that $T'$ really does stabilize to its final value. Thereafter we shall have
  $\Lambda ( e_{0} ,T, \alpha ) >0$. At a later point we shall have all the
  correct trees to apply the Main Lemma once, and these will survive. After
  such a \ point $\Lambda ( e_{0} ,T, \alpha )$ is greater than 1. But only at
  $\beta_{0}$ do we first have $\ensuremath{\operatorname{Liminf}}_{\alpha
  \rightarrow \beta_{0}} \Lambda ( e_{0} ,T, \alpha ) = \omega$ and so
  divergence.
\end{note}

It may already be apparent that the claim that there is an index number $e_{0
}$for the above generalized ittm-recursion can be established readily from the
$\ensuremath{\operatorname{eJ}}$-Recursion Theorem. \ One may argue as
follows, somewhat schematically.

Let $F ( 0,e )$ code the actions of the main programme at (1) above, searching
through increasing $M$-structures for a stable $T'$. (The $e$ is just a dummy
parameter at this stage.) \ With $T'$ stabilized, the programme asks the
oracle query $? \ensuremath{\operatorname{eJ}} ( i,x ) =1/0?$ about $x= \pa{1,
\pa{B_{n}}_{n} ,T,T' ,M}$ with $i=F ( 1,e )$ to be defined next.

Let $F ( k+1,e )$ be the function that returns the index code of the following
blocks of computations:

(1) The actions to do to fulfill the query $Q^{k+1}$, as an explicit
computation. As indicated above these can be listed effectively and the code
of their formal instructions can be given as a function of $k$ - $ q ( k+1 )$
say. This includes the actions to compute the increasing structures and what
to do if stability of any tree passed down subsequently fails. Also included
are, if a stability point is reached that requires a new query to a lower
subcomputation, the actions to collect together the current trees, to form
part of a new coding real $x$.

(2) ``$( x )_{0} := ( x )_{0} +1$'' \ \ [\% Increase the initial index of $x$
by 1 - here to $k+2$.]

(3) The code of the query instruction: ``$? \ensuremath{\operatorname{eJ}} (
\varphi_{e}^{\ensuremath{\operatorname{eJ}}} ( ( x )_{0}  ) ,x ) =0/1?$''.

{\nod}Let the two instructions (2) and (3) have code together $t ( e ) \in
\nat$.

(4) The code of the post-query actions, on receipt of an answer (in the form
of what to do if information is received of a certain kind of tree from a
lower subcomputation becoming unstable etc). \ Again these are effective in
$k$. Let these be $p ( k+1 )$ say.

We thus may loosely represent the total function $F ( k+1,e )$ as:
$$F ( k+1,e ) =
q ( k+1 ) \con t ( e ) \con p ( k+1 ).$$

{\nod}By the $\ensuremath{\operatorname{eJ}}$-Recursion Theorem there is an
index $e_{0}$, so that
$$\varphi^{\ensuremath{\operatorname{eJ}}}_{e_{0}} ( k+1 )  =F ( k+1,e_{0} ) =
q ( k+1 ) \con t ( e_{0} ) \con p ( k+1 ).$$

{\nod}Then our overall computation is: $\{ e_{0}
\}^{\ensuremath{\operatorname{eJ}}} \left( \pa{B_{n} \mid n< \omega} ,T
\right)$.\\

As for the outcome we have as a final claim:

{\nod}{\textbf{Claim:}} {{\em For $A= \bigcup_{n} B_{n} \in \Sigma^{0}_{3}$
and $T$ a recursive subtree of $^{< \omega} \omega$ as above, the programme
$\p{e_{0}} \left( \pa{B_{n} \mid n< \omega} ,T \right)$ will either halt with
a code for a strategy for {\pone}, if such exists, or else will diverge. In
the latter case if it diverges after $\beta$ steps, then a strategy for
{\ptwo} is definable over\/}} $L_{\beta}$.

{\pf} We first observe that the master programme (at $\Lambda =0$) cannot
enter an eventual loop: suppose $( \zeta , \Sigma )$ was its first looping
pair of ordinals. Then the level of computation at times $\zeta$ and $\Sigma$
is the same: $\Lambda ( \zeta ) = \Lambda ( \Sigma ) =0$. But the argument of
Claim 1 of Lemma \ref{4.29}, shows that we must have stability of $T'  $by any
extendible ordinal $\zeta$, and hence, by the specification of $e_{0}$, must
be at a level $>0$ at time $\zeta : \Lambda ( \zeta ) >0$. The same argument
shows that even with $\Lambda ( \zeta ) = \liminf_{\alpha \rightarrow \zeta}
\Lambda ( \alpha ) =0$, we should have $T'$ diminishing unboundedly below the
$2$-extendible $\zeta$ - which cannot happen.

So the computation either halts or diverges. However divergence can only
happen if there is an infinitely descending chain of query calls $Q^{k}$. \
And such has been designed only to happen when we have complete stability of
all our definable trees necessary for the proof of the existence of a
definable winning strategy for {\ptwo} over $L_{\beta}$ - as our procedures
mimic. Lastly the main programme can only halt if it produces a winning
strategy for $\pone .$ \ {\hspace*{\fill}}{\qed Theorem\ref{overall}}\\

Hence by the latter case of the last Claim, strategies for {\ptwo} in such
games are in general not even semi-recursive in
$\ensuremath{\operatorname{eJ}}$.

\begin{corollary}
  There is a procedure $\p{e}$ that only diverges at $\beta_{0}$.
\end{corollary}

{\pf} Let $A= \bigcup_{n< \omega} B_{n} \in \Sigma^{0}_{3}$ be such that $G (
A;T )$ is a win for {\ptwo}, but there is no winning strategy in
$L_{\alpha_{0}}$. Then the computation $\p{e_{0}} \left( \pa{B_{n} \mid n<
\omega} ,T \right)$ above can only diverge at $\beta_{0}$ since a winning
strategy for {\ptwo} is definable over $L_{\beta_{0}}$ but no earlier. \ \
{\hspace*{\fill}}{\qed}

{\bu} An example of such a game, of the type above, is where {\ptwo} must
construct an $\omega$-model of
``$\ensuremath{\operatorname{KP}}+\ensuremath{\operatorname{Det}} (
\Sigma^{0}_{3} )$'', and {\pone} as usual must find a descending chain of
ordinals in {\ptwo}'s model. Then {\ptwo} has an obvious winning strategy, but
there cannot be one where {\ptwo} produces a model with wellfounded part an
ordinal smaller than $\beta_{0}$. We saw in the proof of the theorem above
that the computation in a game of this type, continually constructs codes for
the levels of the $L$-hierarchy unboundedly in $\beta_{0}$, - and hence is
ultimately divergent. We thus have:

\begin{corollary}
  \label{longcomp}There is a program code $f$ so that (i) $\p{f} ( x )$
  computes codes for levels for the $L [ x ]$-hierarchy; (ii) $\p{f} ( 0 )$ is
  divergent , but is not divergent at any stage before $\beta_{0}$, whilst
  computing codes for levels $L_{\alpha}$ for $\alpha$ unbounded in
  $\beta_{0}$.{\hspace*{\fill}} {{\em {\qed}\/}}
\end{corollary}

\begin{corollary}
  $\eta_{0} = \tau_{0}$ - that is Theorem \ref{sigma=tau} holds.
\end{corollary}

{\pf} We have that $\alpha_{0} = \eta_{0}$. By modifying the program of the
last Corollary we can find programs $\p{f} ( 0 )$ which halt cofinally in the
admissible set $L_{\alpha_{0}}$, and hence with ranks of such computations
unbounded in $\alpha_{0}$. Hence $\tau_{0} \geq \alpha_{0}$. By the
Boundedness Lemma \ref{Boundedness} $\tau_{0} \leq \alpha_{0}$.
{\hspace*{\fill}}{\qed}

\begin{lemma}
  Let $a \sset \omega$ be in $L_{\alpha_{0}}$. Then $a$ is
  $\ensuremath{\operatorname{eJ}}$-recursive.{\hspace*{\fill}} {{\em
  {\qed}\/}}
\end{lemma}

The following answers a question of Lubarsky:

\begin{corollary}
  \label{Bob}The reals appearing on the tapes of freezing-ittm-computations of
  {{\em {\cite{Lu10}}\/}} are precisely those of $L_{\beta_{0}}$; similarly
  the supremum of the ranks of the wellfounded parts of divergent computation
  trees is $\beta_{0}$.
\end{corollary}

{\pf}Freezing-ittms computations are, in the terms here, divergent
$\ensuremath{\operatorname{iJ}}$-computations. As
$\ensuremath{\operatorname{eJ}}$ is recursive in
$\ensuremath{\operatorname{iJ}}$ we shall have that the
$\ensuremath{\operatorname{iJ}}$-recursive reals and the
$\ensuremath{\operatorname{eJ}}$-recursive reals coincide. These will be the
reals of $L_{\alpha_{0}}$. By the Boundedness Lemma all such computations are
divergent by $\beta_{0}$, whilst at the same time codes for levels of $L$ for
$\alpha < \beta_{0}$ appear on some $\p{e}$'s tape. Hence the reals appearing
on the divergent $\ensuremath{\operatorname{iJ}}$-computations are those of
$L_{\beta_{0}}$.{\hspace*{\fill}} {\qed}

\begin{corollary}
  The complete semi-decidable-in-$\ensuremath{\operatorname{eJ}}$ set of
  integers
  $$K= \{ ( e,m
  ) \in \omega \times \omega \mid    \ensuremath{\operatorname{eJ}} ( e,m ) =1
  \}$$
  {\nod}is recursively isomorphic to a complete $\game \Sigma^{0}_{3}$ set.
\end{corollary}

{\pf} If $P^{\ensuremath{\operatorname{eJ}}}_{e} ( m )$ is convergent it must
be so before $\beta_{0}$: its convergence is a $\Sigma_{1}$-fact true in
$L_{\beta_{0}}$. By $\Sigma_{1}$-reflection, it is true in $L_{\alpha_{0}}$.
Hence the $\Sigma_{1}$-fact of its convergence is mentioned in the
$\Sigma_{1}$-$\ensuremath{\operatorname{Th}} ( L_{\alpha_{0}} )$. That is $K
\leq_{1} \Sigma_{1}$-$\ensuremath{\operatorname{Th}} ( L_{\alpha_{0}} )
\equiv_{1} S$ where $S$ is a complete $\game \Sigma^{0}_{3}$ set. \ The latter
holds by Theorem \ref{GSigma}. For the converse, we have that $n \in S $if
there is a certain strategy in $L_{\alpha_{0}}$ for a certain game which is
winning for \pone. Such can be found by inspecting the various $L_{\alpha}$
for $\alpha < \alpha_{0}$. And Corollary \ref{longcomp} enables us to run a
computation which is convergent if such can be found. Hence $S \leq_{1} K$.
{\hspace*{\fill}}{\qed}\\

{\nod}{\textbf{Proof of Theorem \ref{2.8}}}

{\nod}The last Corollary proves the (a) (i) iff (iii) direction of the
Theorem, and we have already established (a)(ii) iff (iii) (in the proof of
Theorem \ref{GSigma}). This leaves (b). But this follows from the usual
characterisation of the semi-recursive and co-semi-recursive sets as being
recursive, the admissibility of $L_{\alpha_{0}}$, and that $\alpha_{0} =
\eta_{0}  $.

{\hspace*{\fill}}{\qed} Theorem \ref{2.8}

We may also recast the above arguments as showing:

\begin{corollary}
  Both the theory $T^{1}_{\alpha_{0}}$ and $K$ are
  {\game}$\Sigma^{0}_{3}$-inductive sets of integers.
\end{corollary}

\begin{remark}
  The same considerations show that in fact the whole of
  $\ensuremath{\operatorname{dom}} ( \ensuremath{\operatorname{eJ}} ) \cap
  \omega \times \omega^{< \omega}$ is {\game}$\Sigma^{0}_{3}$-inductive.
\end{remark}

The proofs of Theorems \ref{gamma=beta}, \ref{2.8}, and \ref{sigma=tau} are
now complete (and they cover the statements of the Theorems 1.5-1.8 in Section
1 of the Introduction).

\end{document}